\documentclass[smallcondensed,natbib]{svjour3}
\usepackage{hyperref,amsmath,amssymb,graphicx,color,url,fullpage,mathptmx,subfig,afterpage}

\newcommand{\eps}{\varepsilon}
\newcommand{\f}{\mathbf{f}}
\newcommand{\p}{\mathbf{p}}
\newcommand{\w}{\mathbf{w}}
\renewcommand{\d}{\mathbf{d}}
\renewcommand{\u}{\mathbf{u}}
\renewcommand{\v}{\mathbf{v}}
\renewcommand{\r}{\mathbf{r}}

\renewcommand{\L}{\rm{L}}
\newcommand{\Lib}{\mathcal{L}}
\renewcommand{\H}{\rm{H}}
\newcommand{\V}{\rm{V}}
\newcommand{\A}{\rm{A}}

\title{Boundary-Value Problem Formulations for Computing Invariant Manifolds
and Connecting Orbits in the Circular Restricted Three Body Problem}
\titlerunning{Manifolds and Connecting Orbits in the
Circular Restricted Three Body Problem}
\author{R.~C.~Calleja \and E.~J.~Doedel \and A.~R.~Humphries \and
A.~Lemus-Rodr\'{\i}guez \and B.~E.~Oldeman}
\institute{
R.~C.~Calleja \and A.~R.~Humphries
\at Mathematics and Statistics, McGill University,
805 Sherbrooke O.,
Montr\'eal, Qu\'ebec, H3A 2K6, Canada
\and
E.~J.~Doedel \and B.~E.~Oldeman, \email{boldeman@encs.concordia.ca}
\at Computer Science, Concordia University,
1455 boulevard de Maisonneuve O.,
Montr\'eal, Qu\'ebec, H3G 1M8, Canada
\and
A.~Lemus-Rodr\'{\i}guez
\at Mathematics and Statistics,
Concordia University, 1455 boulevard de Maisonneuve O.,
Montr\'eal, Qu\'ebec, H3G 1M8,
Canada
}
\begin{document}
\bibliographystyle{abbrvnat}
\maketitle
\begin{abstract}
We demonstrate the remarkable effectiveness of boundary value formulations
coupled to numerical continuation for the computation of stable and unstable manifolds
in systems of ordinary differential equations.
Specifically, we consider the Circular Restricted Three-Body Problem (CR3BP), which
models the motion of a satellite in an Earth-Moon-like system. The CR3BP has many
well-known families of periodic orbits, such as the planar Lyapunov orbits and the
non-planar Vertical and Halo orbits.
We compute the unstable manifolds of selected Vertical and Halo orbits, which in
several cases leads to the detection of heteroclinic connections from such a periodic
orbit to invariant tori. Subsequent continuation of these connecting orbits with a
suitable end point condition and allowing the energy level to vary leads to the
further detection of apparent homoclinic connections from the base periodic orbit
to itself, or the detection of heteroclinic connections from the base periodic orbit
to other periodic orbits. Some of these connecting orbits are of potential interest
in space mission design.
\end{abstract}

\keywords{Restricted three body problem \and Boundary value problems
  \and Invariant manifolds \and
  Connecting orbits \and Numerical continuation}
\section{Introduction}
\label{sec:introduction}
Numerical continuation of solutions to boundary value problems (BVPs) has been used
extensively to study periodic orbits and their bifurcations, including homoclinic and
heteroclinic orbits, in a wide variety of systems.
For a recent overview see the articles in~\citet{doedelbday}.
In particular, these techniques have been applied to compute families of periodic orbits
in the Circular Restricted 3-Body Problem (CR3BP); see, for example,~\citet{dpkdgv}.
In recent years invariant manifolds such as those in the Lorenz system have also
been computed in detail using numerical continuation~\citep{krauskopf:07,aguirre:11,doedel:11}.
At the same time, continuation methods have been developed
for computing and continuing homoclinic and heteroclinic connecting orbits between periodic
orbits and equilibria or periodic orbits~\citep{dkkv1,dkkv2,kr}.
In the current work we use a combination of these techniques to illustrate their
effectiveness in computing stable and unstable manifolds and connecting orbits in the CR3BP.

There is much literature on invariant manifolds and connecting orbits
in the CR3BP; see for example
\citet{Gom-Koo-Lo-Mar-Mas-Ros-04,marsden,Lo-Ros-97,
Dav-And-Sche-Bor-10,Dav-And-Sche-Bor-11,Tan-Fan-Ren-Per-Go-Mas-10}.
In particular, connecting orbits in the planar CR3BP are well understood.
The existence of connecting orbits in the planar problem has been proved analytically
in \citet{Lli-Mar-Sim-85}, and by computer assisted methods in \citet{Wil-Zgl-03,Wil-Zgl-05}.
Furthermore, these orbits have been
extensively studied numerically using initial-value
  techniques and semi-analytical tools; see
\citet{Bar-Mon-Oll-09,Can-Mas-06} and references therein.
In the case of initial-value techniques the initial conditions are varied in order
for an appropriately chosen end point condition to be satisfied. This approach is
commonly referred to as a ``shooting method'' and, for a more stable version,
``multiple shooting''.
Initial-value techniques can also be very effective in the computation of invariant 
manifolds in the CR3BP. However, sensitive dependence on initial conditions may
leave parts of the manifolds unexplored, unless very high accuracy is used.
Other efficient methods for computing invariant manifolds include semi-analytical
approximations \citep{Jor-Mas-99,Ale-Gom-Mas-09,Gom-Mon-01}.
The latter methods are very precise in a neighborhood of the center of expansion,
and rely on other methods to extend the manifolds outside these neighborhoods
\citep{Gom-Jor-Mas-Sim-01}.

Invariant manifold techniques around libration points have been used successfully
in mission design \citep{Lo-And-Whi-Rom-04}.
The Genesis spacecraft mission, designed to collect samples of solar wind and return
them to the Earth
 \citep{Lo-Wil-Bol-Han-Hah-Hir-01}, is often considered as the
first mission to use invariant manifolds
for its planning, while other missions have used libration point techniques
\citep{Dun-Far-03}. Having a precise idea of the geometry of invariant manifolds and their connections
is desirable in the design of complex low thrust missions.

Using our continuation approach we construct a continuous solution family
in the manifold as the
initial value is allowed to vary along a given curve.
The continuation step size governs the distance between any computed trajectory
and the next trajectory to be computed. Here ``distance'' includes the change in
the entire trajectory, and not only in the initial conditions. In fact, this
distance typically also includes other variables in the continuation process,
such as the integration time and the arclength of the trajectory. This formulation
allows the entire manifold to be covered, up to a prescribed length, integration
time, or other termination criterion, even in very sensitive cases.
Special orbits, such as connecting orbits to saddle-type objects, can be detected
during the continuation.
For example, a straightforward computation of the Lorenz manifold in this fashion
yielded up to 512 connecting orbits having extremely close initial values, as
reported in \citet{doedel:06}. In related work, the intersections of the Lorenz
manifold with a sphere are studied in \citet{aguirre:11} and \citet{doedel:11}.
In the case of the Lorenz manifold, the sensitivity on initial conditions results
from the significant difference in magnitude of the two real negative eigenvalues
of the zero equilibrium that give rise to this manifold. Fixed precision
integration in negative time may cover only part of this manifold, missing a portion
in and near the direction of the eigenvector of the smaller negative eigenvalue.

In this paper we give an overview of continuation techniques, as used to compute
periodic orbits, invariant manifolds, and connecting orbits. We also give several
examples that illustrate how these techniques provide an effective and relatively
easy-to-use tool for exploring selected portions of phase space.
The richness of the solution structure of the CR3BP limits the extent of our
illustrations. However, the techniques presented here are expected to be useful
in further studies. In this respect, the current version of the freely available
AUTO software \citep{AUTO} includes demos that can be used to re-compute some of
the numerical results presented here, including their graphical representation.
These demos can also be adapted relatively easily to perform similar numerical
studies of stable and unstable manifolds of other periodic orbits in the CR3BP that
are not considered in this paper, as well as for entirely different applications.

This paper is organized as follows. In Sect.~\ref{sec:cr3bp} we recall some
well-known facts about the CR3BP, namely, its equilibria, the libration points,
and the basic periodic solution families that will be considered in this paper, namely
the planar Lyapunov orbits and the Vertical, Halo and Axial families.
In Sect.~\ref{sec:periodic-orbits} we review how boundary value techniques are used
to compute periodic orbits in conservative systems, and how these techniques can also
be used to compute the eigenfunctions associated with selected Floquet multipliers.
These data provide a linear approximation of the unstable manifold of the periodic orbit.

Sect.~\ref{sec:manifold-algorithm} describes the continuation method used for 
computing unstable manifolds of periodic orbits. This involves first setting
up an extended system with both the periodic orbit and its eigenfunction.
Using the resulting information an initial orbit within the manifold is
computed. This orbit is then continued, as its starting point is free to vary
along a line that is tangent to a linear approximation of the unstable manifold,
thereby tracing out the manifold.
The algorithm, using pseudo-arclength continuation,
is not guaranteed to compute the whole manifold in a single computation,
because obstacles may be encountered during the continuation. In such cases 
the manifold may be completed, for example, by additional continuations from 
different starting orbits, or by using a suitably adapted continuation procedure 
as is done in \citet{doedel:11}.  However, the obstacles themselves are
also of interest. They may correspond to orbits in the unstable manifold
that require an arbitrary long time interval to reach a specified termination
plane, because they pass arbitrarily close to a connecting orbit between the
original unstable periodic orbit and another invariant object. In this way 
connecting orbits can be detected, as was done to detect the 512 heteroclinic 
connections presented in \citet{doedel:06}.

In Sect.~\ref{sec:manifold-examples} we show the results of computations of
unstable manifolds of Vertical $\V_1$ and Halo $\H_1$ orbits. When the manifolds
are computed to a sufficient distance from the original periodic orbit, we find
what appear to be heteroclinic connecting orbits from the original periodic orbit
to an invariant torus. The tori found this way must have saddle-type instability,
since the connecting orbit approaches it, but ultimately also leaves the
neighborhood of the torus. Such connecting orbits may be more difficult to find
with initial value integration.

In Sect.~\ref{sec:connecting-orbits} we describe a method for continuing the
periodic-orbit-to-torus connections as solutions of an $18$-dimensional ODE,
when the energy is allowed to vary. These computations lead to the detection
of other interesting connecting orbits. Sect.~\ref{sec:connecting-orbits-examples}
discusses three representative examples.
First we consider a family of connections from $\H_1$ Halo orbits. These
connections loop once around the Earth before approaching an invariant torus,
which is itself close to the $\H_1$ orbit. We refer to such a torus as a
quasi-$\H_1$ torus. During the continuation, with changing energy of the
originating $\H_1$ Halo orbit, we encounter a number of interesting connecting
orbits. Specifically, we find a homoclinic orbit from an $\H_1$ Halo orbit
to itself that loops once around the Earth, a heteroclinic connection from a
northern $\H_1$ Halo orbit to its southern counterpart, a heteroclinic connection
to a planar $\L_1$ Lyapunov orbit, and a connection to a 5:1 resonant orbit on
a torus near the corresponding southern Halo orbit. We also find a connecting orbit
from an $\H_1$ Halo orbit to a torus on which the orbit bounces back and forth
between a northern Axial $\A_1$ orbit and its southern counterpart. Each of these
special connecting orbits occurs for a specific energy of the originating $\H_1$
Halo orbit (and of course, by conservation of energy, the orbit it connects to
has the same energy). Secondly we study a family of connecting orbits from an 
$\H_1$ Halo orbit that loop four times around the Earth. We find connections to
an $\L_2$ Lyapunov orbit, an $\H_2$ Halo orbit, and to a 5:1 and a 6:1 resonant
torus near the libration point $\Lib_2$. Thirdly, we consider a family of
connecting orbits on the Moon-side of the unstable manifold of the $\H_1$ Halo
orbits that connect directly (without looping around the Earth) to a torus near
$\Lib_2$. We find an example of a direct connecting orbit from an $\H_1$ Halo
orbit to a planar $\L_2$ Lyapunov orbit. To the best of our knowledge, these
connecting orbits have not been found before. We must stress, however, that from
a space mission design point of view, these orbits are sensitive to initial
conditions and would require control techniques to stay on them.

In Sect.~\ref{sec:discussion}, we discuss global theoretical aspects of our results 
and their relation to the existing literature.  In particular we see that
the connecting orbits from $\H_1$ Halo orbits to $\A_1$ Axial orbits
or to $\L_1$ or $\L_2$ planar Lyapunov orbits are codimension-one in
the dynamical systems sense, and hence should occur for specific values of the energy
of the originating $\H_1$ Halo orbit, as we observed numerically. In
contrast, homoclinic and heteroclinic connecting orbits between $\H_1$
Halo orbits, which were observed numerically, are codimension-two, and so should
not normally occur. We show that the connecting orbits from the $\H_1$
Halo orbits to quasi-$\H_1$ tori are generic, and suggest that the
numerically observed homoclinic and heteroclinic connecting orbits
between $\H_1$ Halo orbits are actually connections to quasi-$\H_1$ tori where the minor
radius of the torus is so small as to make the torus visibly (and for the purpose of
space mission design) indistinguishable from the $\H_1$ Halo orbit that it envelopes.

Finally in Sect.~\ref{sec:numerical-aspects} we discuss some computational
aspects of our numerical computations, such as the discretization used, and
typical computer time needed.
\section{The Circular Restricted 3-Body Problem}
\label{sec:cr3bp}

\begin{figure}[htb]
\begin{center}
{\setlength{\unitlength}{3947sp}\fontsize{12}{14.4pt}
\subfloat[]{
\begin{picture}(0,0)
\includegraphics{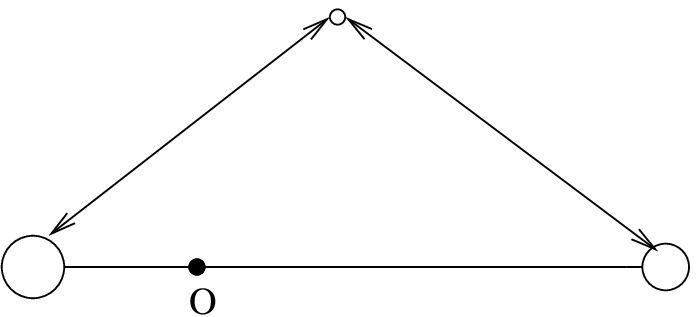}
\end{picture}
\begin{picture}(3315,1662)(2640,-5401)
\put(3151,-4600){$r_1$}
\put(5101,-4561){$r_2$}
\put(4276,-5311){$1-\mu$}
\put(3076,-5311){$\mu$}
\put(2701,-5236){$E$}
\put(5740,-5220){\fontsize{10}{12.0pt}$M$}
\put(4276,-3886){$S$}
\end{picture}}
\hspace{1cm}
\subfloat[]{
\begin{picture}(0,0)
\includegraphics{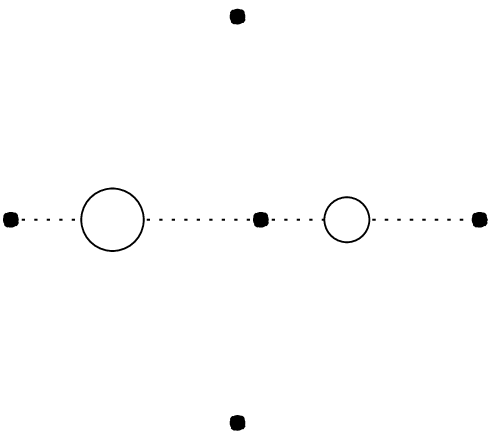}
\end{picture}
\begin{picture}(2355,2145)(1890,-7621)
\put(4150,-6736){$\Lib_2$}
\put(1876,-6736){$\Lib_3$}
\put(3076,-7561){$\Lib_5$}
\put(3076,-5611){$\Lib_4$}
\put(3100,-6736){$\Lib_1$}
\put(2326,-6586){$E$}
\put(3457,-6565){\fontsize{10}{12.0pt}{$M$}}
\end{picture}}}
\caption{(a): Schematic representation of the Circular Restricted Three-Body Problem.
(b): The five libration points.}
\label{fig:r3bp}
\end{center}
\end{figure}

The CR3BP describes the motion of a satellite $S$ with negligible mass
in three-dimensional physical space. The motion is governed by the
gravitational attraction of two heavy bodies, which are assumed to rotate
in circles around their common center of mass; see Fig.~\ref{fig:r3bp}(a).
In this paper we call the heaviest body the ``Earth'' $E$, and the other
heavy body the ``Moon'' $M$, and we use their actual mass ratio, namely,
$\mu\approx 0.01215$. Any other mass ratio is allowed, such as for
the Sun-Jupiter system, with a mass ratio of $\mu\approx 0.0009537$.
Without loss of generality the total mass can be scaled to $1$, so that
the Earth and Moon have mass $0.98785$ and $0.01215$, respectively.

The equations of motion of the CR3BP as given in~\citet{Danby92} are
\begin{equation}
\begin{split}
\ddot x&~=~ 2
\dot y ~+~ x ~-~(1-\mu)~\frac{x+\mu}{{r_1}^3} ~-~ \mu~\frac{x-1+\mu}{{r_2}^3}~,\\
\ddot y &~=~ -2\dot x ~+~ y ~-~ (1-\mu)~\frac{y}{{r_1}^3} ~-~ \mu~\frac{y}{{r_2}^3}~, \\
\ddot z &~=~ -(1-\mu)~\frac{z}{{r_1}^3} ~-~ \mu~\frac{z}{{r_2}^3}~,
\end{split}
\label{eq:cr3bp}
\end{equation}
where $(x,y,z)$ is the position of the zero-mass body, and where
\[r_1=\sqrt{(x+\mu)^2+y^2+z^2}~,\quad r_2=\sqrt{(x-1+\mu)^2+y^2+z^2}~,\]
denote the distance from $S$ to the Earth and to the Moon, respectively.
The CR3BP has one integral of motion, namely, the energy $E$:
$$ E ~=~ \frac{\dot x^2 + \dot y^2+\dot z^2}{2} ~+~ U(x,y,z)~, \qquad
U(x,y,z) ~=~ -\frac{1}{2} (x^2+y^2) ~-~ \frac{1-\mu}{r_1} ~-~ \frac{\mu}{r_2}
	~-~ \mu \frac{1-\mu}{2}~,$$
where $U(x,y,z)$ is the effective potential.
Astronomers also often use the \emph{Jacobi constant}~ C, defined as $C=-2E$.

\begin{figure}[htbp]
\begin{center}
\subfloat[]{
\includegraphics{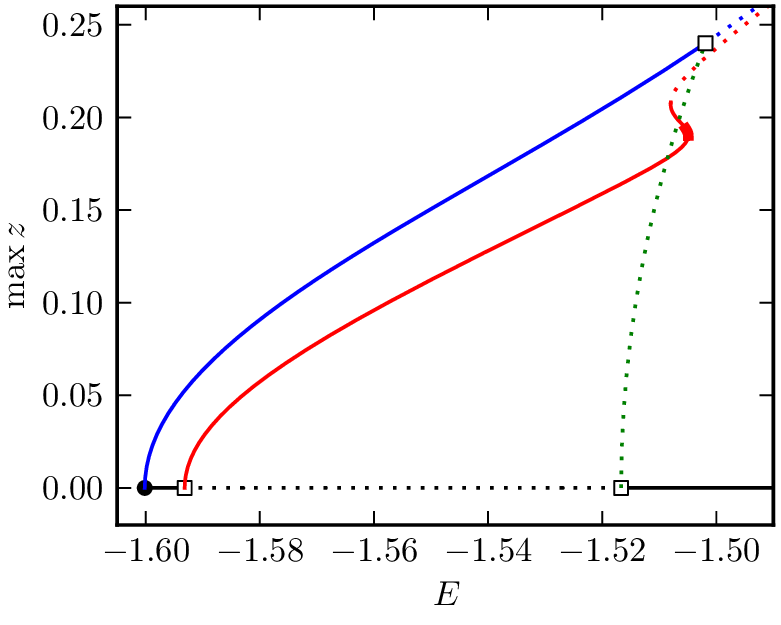}
\put(-186,32){$\Lib_1$}
\put(-126,44){$\L_1$}
\put(-126,77){$\H_1$}
\put(-126,113){$\V_1$}
\put(-41,77){$\A_1$}
\put(-18,160){$\V_{11}$}
\put(-169,32){$\L_{11}$}
\put(-43,32){$\L_{12}$}
}
\subfloat[]{
\includegraphics{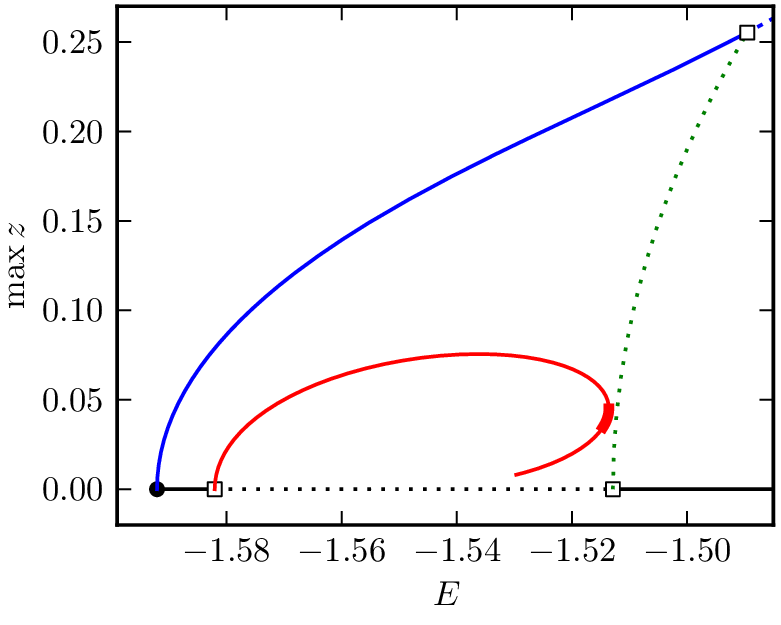}
\put(-183,32){$\Lib_2$}
\put(-136,42){$\L_2$}
\put(-136,74){$\H_2$}
\put(-136,113){$\V_2$}
\put(-36,100){$\A_2$}
\put(-25.5,171){$\V_{21}$}
\put(-160,32){$\L_{21}$}
\put(-46,32){$\L_{22}$}
}
\caption{Bifurcation diagrams of families of periodic solutions of the
Earth-Moon system bifurcating from the libration points (a) $\Lib_1$ and (b) $\Lib_2$.
A detailed description of the families represented in this diagram can be found in
\citet{drpkdgv}. The ``thick'' portions of the curves labeled $\H_1$ and $\H_2$
denote periodic orbits where all 6 Floquet multipliers are on the unit
circle: $1, 1, \exp(\pm i c), \exp(\pm i d)$.
For thin solid portions exactly two real Floquet multipliers are off the
unit circle: $1, 1, a, 1/a, \exp(\pm i c)$.
The dotted portions denote periodic orbits where all
6 Floquet multipliers are real: $1, 1, a, 1/a, b, 1/b$.
Here $a, b, c, d\in \mathbb{R}$, $|a|>1, |b|>1$, and $c, d \in (0,\pi)$.
The small squares labeled $\L_{ij}$ and $\V_{i1}$ denote branch points.}
\label{fig:families}
\end{center}
\end{figure}
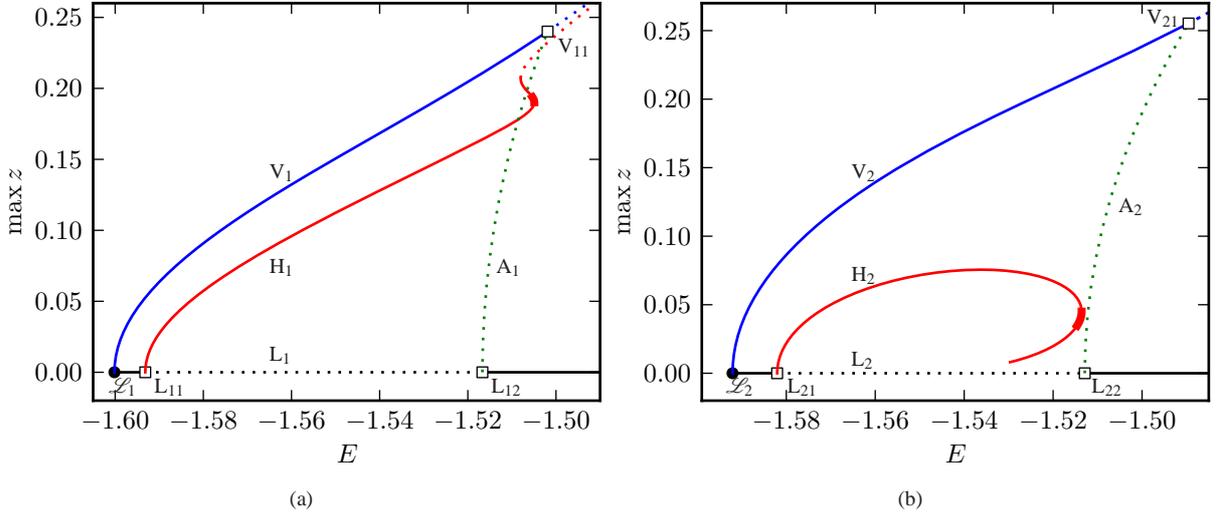
\begin{figure}[htbp]
\begin{center}
\subfloat[]{
\includegraphics[scale=0.29,clip,trim=0mm 5mm 0mm 20mm]{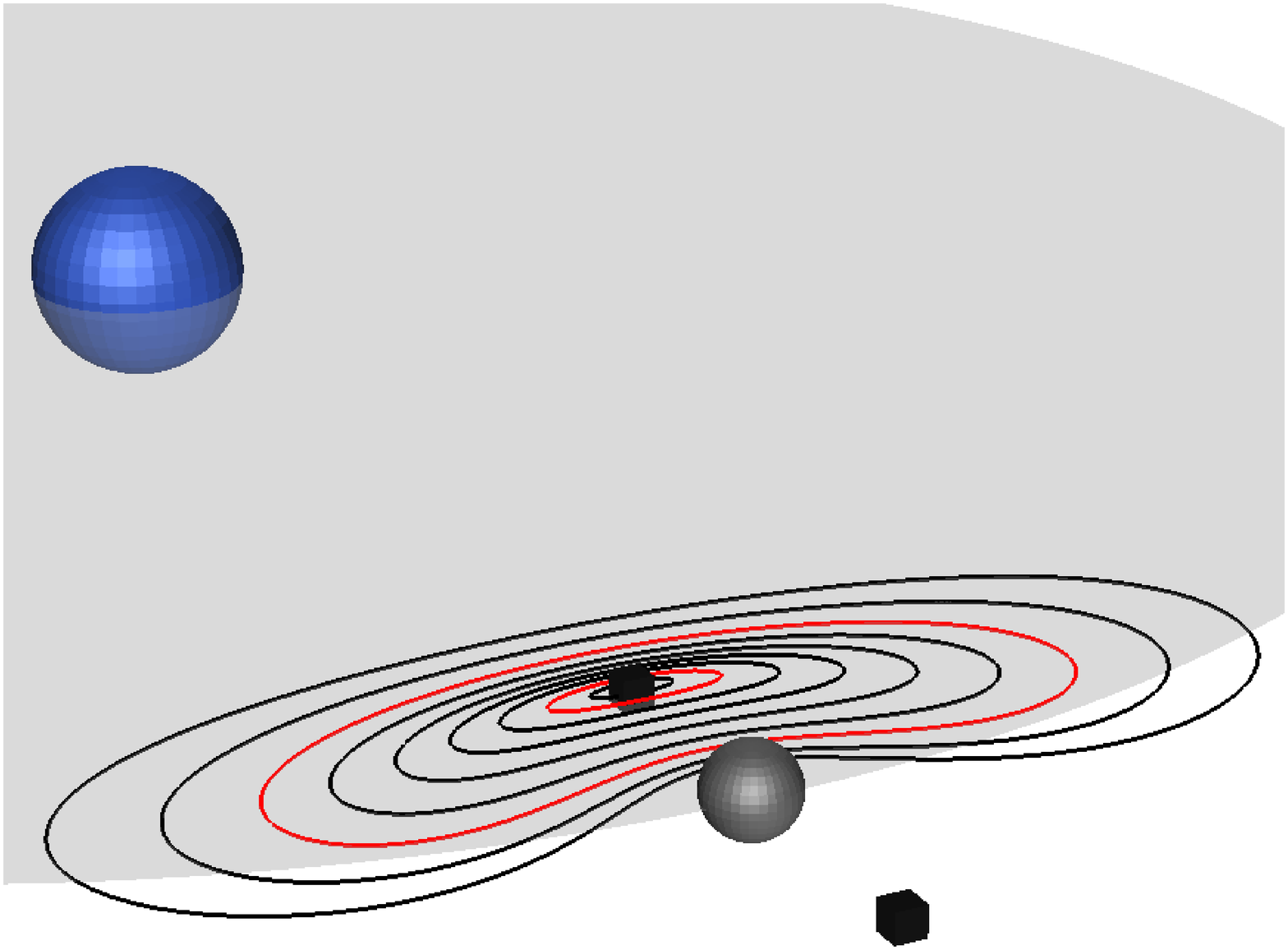}
\put(-143,69){\vector(0,1){40}}
\put(-147,63){$\L_{11}$}
\put(-180,53){\vector(0,1){30}}
\put(-184,47){$\L_{12}$}
}
\subfloat[]{
\includegraphics[scale=0.29,clip,trim=40mm 5mm 0mm 20mm]{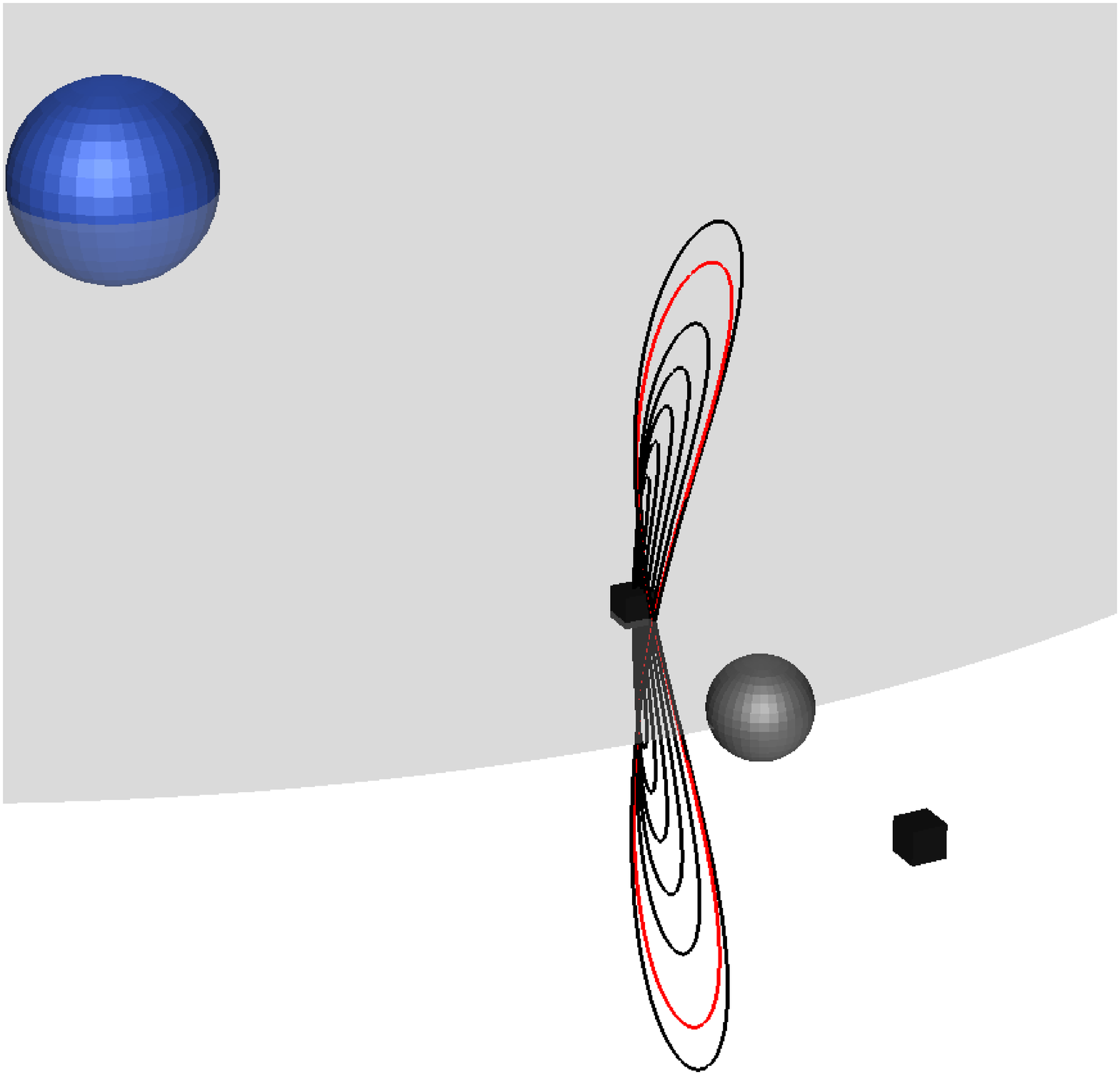}
\put(-108,6){\vector(1,1){10}}
\put(-117,1){$\V_{11}$}
}
\end{center}
\caption{
(a):
A selection of periodic orbits from
the planar Lyapunov family $\L_1$, which bifurcates from the libration
point $\Lib_1$, as seen in the bifurcation diagram in
Fig.~\ref{fig:families}(a).
Labeled are the special Lyapunov orbits $\L_{11}$ from which the Halo family $\H_1$ emanates,
and $\L_{12}$ from which the Axial family $\A_1$ bifurcates.
In this and subsequent figures the small cubes denote libration points.
(b): Selected orbits from the Vertical family $\V_1$, which
also bifurcates from the libration point $\Lib_1$.
In the linear approximation near $\Lib_1$
these orbits are indeed ``vertical'', {\it i.e.}, $x$ and $y$ are constant along
it, with $y=0$, while $z$ oscillates around zero. The Axial family bifurcates from
the orbit labeled $\V_{11}$.}
\label{fig:l1v1}
\end{figure}
\begin{figure}[htbp]
\begin{center}
\subfloat[]{
\includegraphics[scale=0.29,clip,trim=40mm 25mm 0mm 0mm]{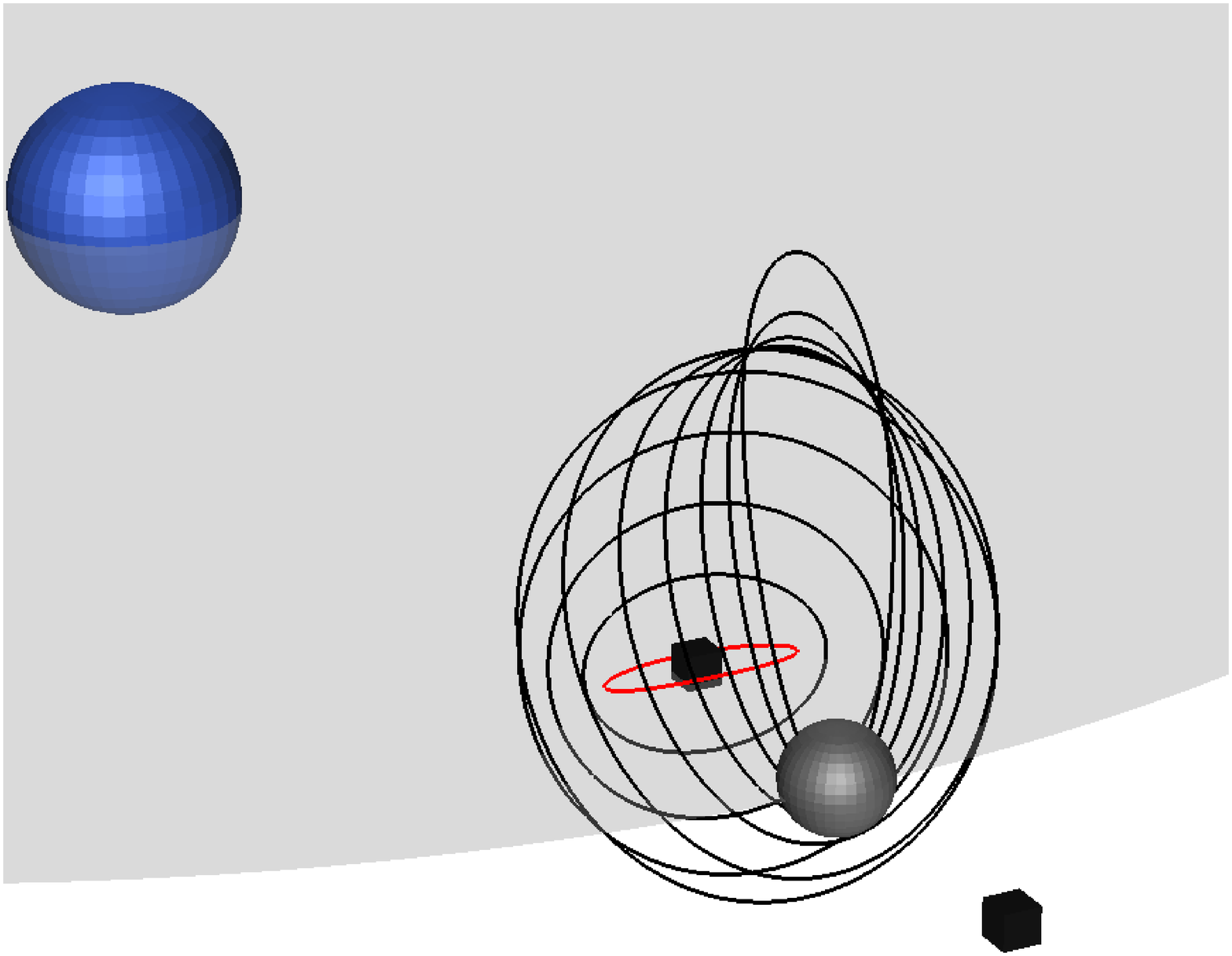}
\put(-125,40){\vector(0,1){50}}
\put(-129,34){$\L_{11}$}
}
\subfloat[]{
\includegraphics[scale=0.29,clip,trim=40mm 5mm 0mm 20mm]{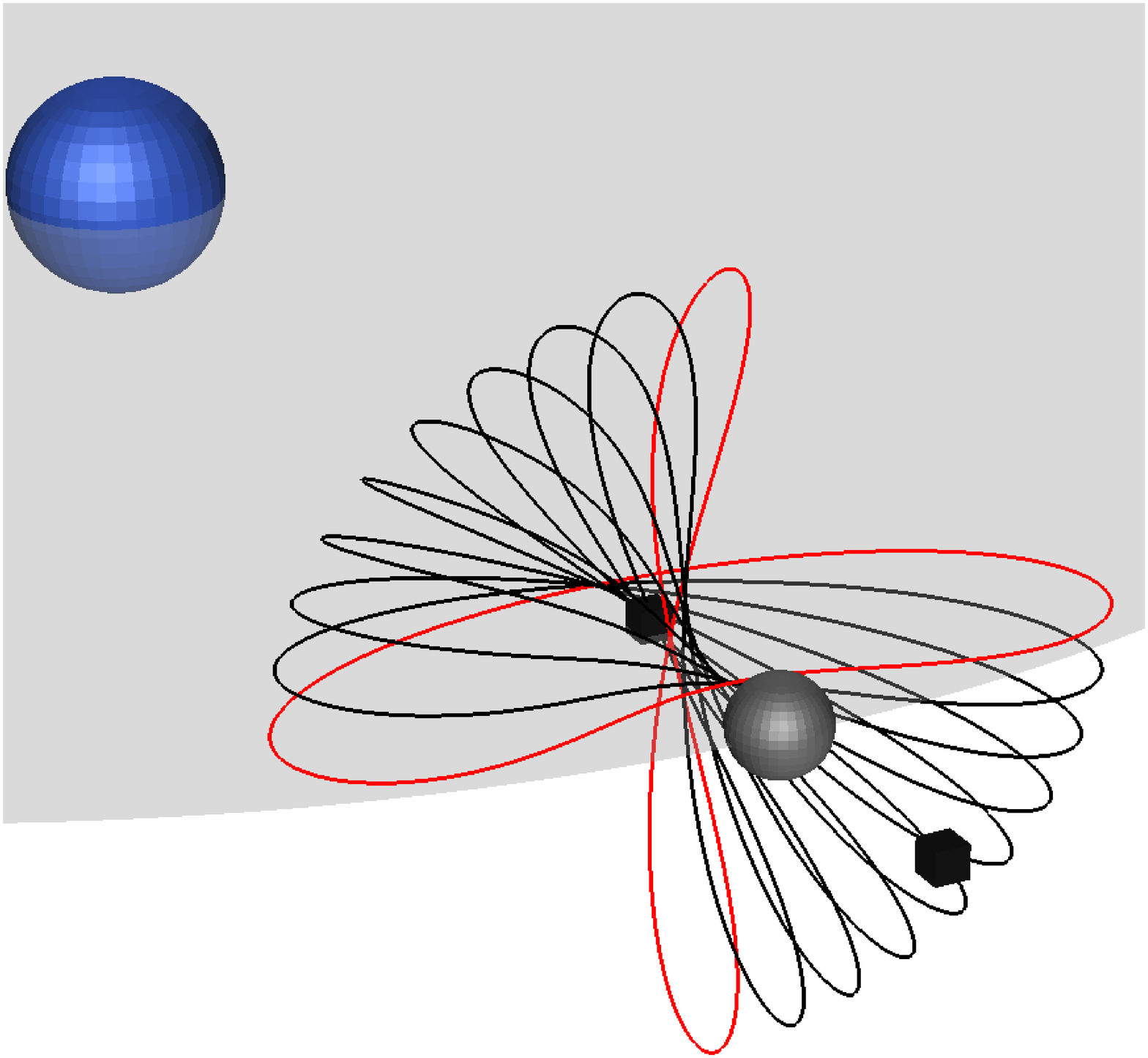}
\put(-138,73){$\L_{12}$}
\put(-118,25){$\V_{11}$}
}
\caption{
(a): Selected orbits from the Halo family $\H_1$ that bifurcates from the Lyapunov
family $\L_1$ at $\L_{11}$.
(b): Selected orbits from the Axial family $\A_1$ that connects to the Lyapunov family
$\L_1$ at $\L_{12}$ and to the Vertical family $\V_1$ at $\V_{11}$: see also Fig.~\ref{fig:families}.}
\label{fig:h1a1}
\end{center}
\end{figure}

%
Libration points, in both the planar and three-dimensional spatial system,
are equilibrium points in a co-rotating frame; see Fig.~\ref{fig:r3bp}(b).
As already used, we denote the libration points by $\Lib_1$, $\Lib_2$,
~$\cdots$~, $\Lib_5$.
There are families of periodic orbits (in the co-rotating frame) that bifurcate
from each of these libration points, and we refer to these as the {\it primary
families}.
Many more families subsequently bifurcate from the primary families.
We refer to the bifurcation points as {\em branch points}.
Several families of periodic solutions of the Earth-Moon system are represented in
Fig.~\ref{fig:families}; see also \citet{dpkdgv,drpkdgv} and
references therein.
In the present work we focus on four families, namely,
the Vertical orbits $\V_i$, the planar Lyapunov orbits $\L_i$,
the Halo orbits $\H_i$, and the Axial orbits $\A_i$.
The families of planar Lyapunov orbits $\L_i$ and the Vertical orbits
$\V_i$ emanate directly from the libration points $\Lib_i$; these are primary
families (Fig.~\ref{fig:l1v1}).
The families of Halo orbits bifurcate from the families of Lyapunov orbits $\L_i$,
while the families of Axial orbits connect and bifurcate from the $\V_i$ and $\L_i$ families.
We refer to the Halo and Axial orbits as {\it secondary families}
(Fig.~\ref{fig:h1a1}).
These families are all well-documented in the literature, but their
names are sometimes different. For example,
the Halo, Axial, and Vertical orbits are known as type ``A'', ``B'', and
``C'', respectively, in \citet{Goudas-61} and \citet{Henon-73}.
\citet{Farquhar-68} coined the name ``Halo'' for that family.
The term ``Axial'' comes from \citet{drpkdgv}, whereas \citet{dpkdgv}
used the term ``Y'' for ``Yellow''.

%
%
\section{BVPs for Periodic Orbits and Eigenfunctions}
\label{sec:periodic-orbits}

In this and the next section, we describe the algorithms used in the various
stages of the computations in some detail, so that it will be possible for the
reader to replicate the algorithms in other applications, and to better understand
the functioning of the downloadable CR3BP demos \citep{AUTO}.
Specifically, in this section we explain the preliminary computations that precede
the actual computation of stable and unstable manifolds of periodic orbits, namely,
the computation of the periodic orbits themselves and of their associated
eigenfunctions.  The discussion follows that in \citet{dpkdgv,dkkv1}.

To formulate a suitable boundary value problem in AUTO, the second order system
of ODEs \eqref{eq:cr3bp} first needs to be rewritten as a six-dimensional first
order system,
$$ \dot{\u}(t) ~=~ \hat\f(\u(t), \mu)~, \qquad \hat\f:\mathbb{R}^6 \times \mathbb{R}
	\to \mathbb{R}^6~,$$
where $\u=(x,y,z,v_x,v_y,v_z)=(x,y,z,\dot{x},\dot{y},\dot{z})$.
As usual, when we plot the orbits graphically we project into $\mathbb{R}^3$
by only plotting the spatial coordinates $(x,y,z)$.
Time is scaled to the interval $[0,1]$ in the BVP formulation for computing a
periodic orbit, which changes the system of differential equations to
$$ \dot\u(t) ~=~ T~\hat\f(\u(t), \mu)~,$$
where $T$ is the period of the periodic orbit.
In addition, for a conservative system with one conserved quantity, we need to
add a term with an  ``unfolding parameter'' in order for the BVP continuation
computations to be formally well-posed \citep{drpkdgv,afgdv};
see also \citet{mfgv07}.
A suitable and convenient choice in the specific case of the CR3BP is the term
$\sigma \d(\u)$, where $\d(\u) = (0,0,0,v_x,v_y,v_z)$.
The vector field with the unfolding term then becomes
$$ \dot{\u}(t) ~=~ T~\hat\f(\u(t), \mu) ~+~ \sigma \d(\u(t))~,$$
which from here on we simply write as
\begin{equation}
\dot{\u}(t) ~=~ T~\f(\u(t),\sigma)~,
\label{eq:ODE}
\end{equation}
also omitting the mass-ratio parameter $\mu$, as it is typically fixed in
the computation of families of periodic orbits.
Notice that the specific choice of unfolding term $\d(\u)$
used here would represent a damping (or forcing) term if $\sigma\ne0$, which would
preclude the existence of periodic orbits. However,
the unfolding parameter $\sigma$ is one of the unknowns in the continuation
procedure, and will always be zero (to numerical precision) once solved for.
Thus the unfolding term is simply a technical device necessary to obtain well-posedness
of the BVP, and we do not force or damp the equations of motion.

Written in full, the system is therefore given by
\begin{equation}
\begin{split}
\dot x &= T~v_x~,\\
\dot y &= T~v_y~,\\
\dot z &= T~v_z~,\\
\dot v_x &= T~[2 v_y + x - (1-\mu) (x+\mu) r_1^{-3}
              - \mu (x-1+\mu) r_2^{-3}] + \sigma v_x~,\\
\dot v_y &= T~[-2 v_x + y - (1-\mu) y r_1^{-3}
              - \mu y r_2^{-3}] + \sigma v_y~,\\
\dot v_z &= T~[-(1-\mu)z{{r_1}^{-3}} - \mu z{{r_2}^{-3}}] + \sigma v_z~.
\end{split}
\label{eq:writtenout}
\end{equation}
To complete the BVP formulation we need to add the periodicity equations
\begin{equation}
\u(1) ~=~ \u(0)~.
\label{eq:bc}
\end{equation}
If $\u(t)$ solves Eqs.~\eqref{eq:writtenout} and~\eqref{eq:bc}
then $\u_\alpha(t)=\u(t+\alpha)$ is also a solution for any 
time-translation $\alpha$.
To specify a unique solution we impose the
phase constraint \citep{AUTO81}
\begin{equation}
\int_0^1 \langle \u(t),\dot{\u}_0(t) \rangle dt ~=~0~,
\label{eq:phase}
\end{equation}
where $\u_0(t)$ is the preceding computed solution along the solution family.
Furthermore, for the purpose of continuing a family of periodic solutions,
we add Keller's pseudo-arclength constraint \citep{hbk:77}, which in the
current setting takes the form
\begin{equation}
\int_0^1
       \langle \u(t)-\u_0 (t), \u_0' (t)\rangle dt
       + (T-T_0) T_0'
       + (\sigma -\sigma_0) {\sigma}_0'
       ~=~ \Delta s~,
\label{eq:pseudo}
\end{equation}
where $(\u_0, T_0, \sigma_0)$ corresponds to a computed solution
along a solution family, $(\u, T, \sigma)$ is the next solution to be
computed, and $\Delta s$ is the continuation step size. The notation
``$~'~$'' denotes the derivative with respect to $\Delta s$ at $\Delta s=0$
and $\langle \cdot,\cdot \rangle$ denotes the dot product.
Since we are dealing with a conservative system, we already have families of
periodic solutions even when the mass-ratio $\mu$ is fixed.
For a computed solution $(\u_0, T_0, \sigma_0)$, and a given step $\Delta s$,
the unknowns to be solved for in any continuation step are the periodic
orbit $\u(t)$, its period $T$, and the unfolding parameter $\sigma$.
Eq.~\ref{eq:pseudo} forces all of these unknowns to be close to
those of the previous solution. In particular, $\u(t)$ must be close to
$\u_0(t)$ for \emph{all} $t$, and not just for $t=0$.
We reiterate that the unfolding parameter $\sigma$ is an active unknown in
the computations that regularizes the boundary value formulation, although
it will be found to be zero to numerical precision. During the
continuation the Floquet multipliers of the periodic solution are
monitored by computing a special decomposition of the monodromy matrix
that arises as a by-product of the decomposition
of the Jacobian of the collocation system \citep{FaJe:91}.

For the purpose of computing a stable or unstable manifold of a periodic
orbit we need the corresponding Floquet eigenfunction. Specifically,
we assume that the periodic orbit has a single, real, positive Floquet
multiplier outside the unit circle in the complex plane, which gives
rise to a two-dimensional unstable manifold of the periodic orbit in phase
space.
The eigenfunction corresponding to this multiplier provides a linear
approximation to the manifold close to the periodic orbit.
In \citet{dkkv1} it is shown that this eigenfunction can be obtained as a
solution $\v(t)$ of the BVP
\begin{equation}
\begin{split} 
&\dot\v(t) = T \f_{\u} \bigl(\u(t),0\bigr) \v(t)+\lambda \v(t)~,\\
&\v(1) = \pm\v(0)~, \\
&\langle \v(0),\v(0) \rangle =\rho~,
\end{split}
\label{eq:eigenfunction}
\end{equation}
where $\v(1) = +\v(0)$ in the case of a positive multiplier,
and $\v(1) =-\v(0)$ in the case of a negative multiplier.
Here, $\lambda$ is the characteristic exponent and
the corresponding Floquet multiplier is given by $\pm e^\lambda$.
Eq.~\eqref{eq:eigenfunction} represents the Floquet eigenfunction/eigenvalue
relation for the linearization of Eq.~\eqref{eq:ODE} about a periodic
orbit $\u(t)$.
The norm of the value of the eigenfunction at time $t=0$ is normalized
to be $\sqrt{\rho}$, where typically we use $\rho=1$.
If only one Floquet multiplier is real and greater than one in absolute value,
this gives a unique (up to sign) unstable eigenfunction $\v(t)$.
Likewise, a unique stable eigenfunction is obtained if only one Floquet
multiplier is real and less than one in absolute value.
In our illustrations we only compute unstable manifolds, so we have
$\lambda>0$ in Eq.~\eqref{eq:eigenfunction}; however all algorithms apply equally
well to stable manifolds. We also restrict to the case
where the Floquet multiplier of interest is positive,
so $\v(1) = +\v(0)$ in Eq.~\eqref{eq:eigenfunction},
and the corresponding manifold is orientable rather than twisted.
A linear approximation of the unstable manifold at time zero is then
given by
\begin{equation}
  \r(0) ~=~ \u(0)+\eps \v(0)~,
\end{equation}
for $\eps$ small.

An alternative formulation is to put the actual Floquet multiplier
in the boundary condition rather than in the linearized differential
equation, using the variational equation
$\dot\v(t) = T\f_{\u}(\u(t),0) \v(t)$
and boundary condition
$\v(1) = \kappa\v(0)$, where $\kappa$ is the actual Floquet multiplier.
However, the formulation in Eq.~\eqref{eq:eigenfunction}, as used here, has
been found to be more appropriate
for numerical purposes \citep{dkkv1}. This is related to the fact that the multipliers,
{\it i.e.}, the values of $\kappa=e^\lambda$, can be very large or very small.

\noindent The algorithmic steps that lead to the linear approximation of the
unstable manifold are then as follows; here described for the case
of a Halo orbit in the $\H_1$ family.
\begin{enumerate}
\item
The libration points, which are the equilibria of Eq.~\eqref{eq:cr3bp},
are easily determined \citep{Szebehely67}. In our continuation context
we note that they have zero velocity components, as well as $z=0$, and
that for varying $\mu$ their $x$ and $y$ components lie on connected
curves, as shown in Fig.~\ref{fig:libration}.
Starting from, for example, the curve of equilibria
$x^2+y^2=1$, that exists when $\mu=0$ (the curve containing the point $q$ in
Fig.~\ref{fig:libration}), the libration points bifurcate from $x=1/2$, $y=\pm\sqrt{3}/2$
and $y=0$, $x=\pm1$ and we can reach each of the libration points at any
given nonzero value of $\mu$ via a connected path. 
The eigenvalues of the target libration point(s) are also computed.
\begin{figure}[htb]
\begin{center}
\includegraphics{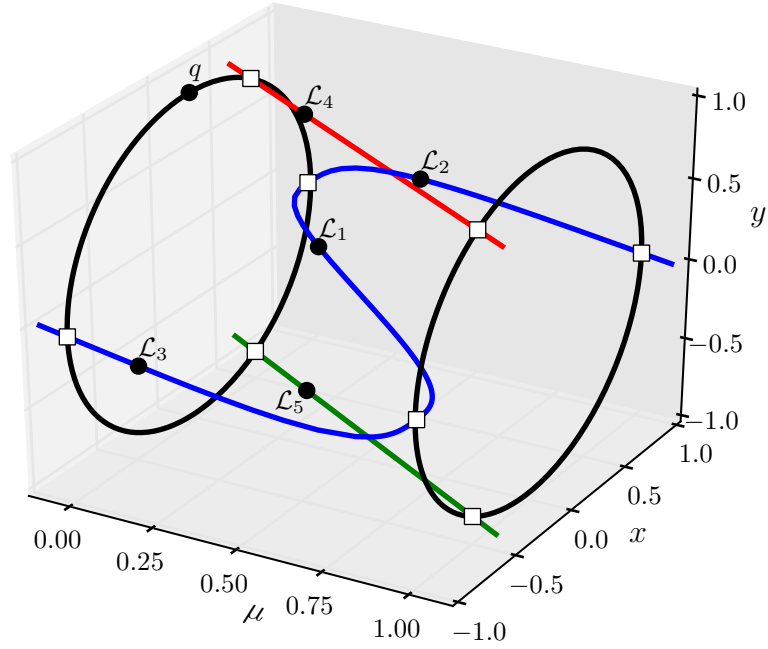}
\end{center}
\caption{
Computation of the libration points. The point $q$ on the circle of
equilibria is a suitable starting point for determining the libration points for
a given nonzero value of $\mu$ by continuation. For clarity, the five libration 
points are indicated in the diagram for the case $\mu=0.25$, rather than for the 
Earth-Moon case ($\mu=0.01215$), as $\mu$ is close to zero in the latter case.
The small squares are branch points.}
\label{fig:libration}
\end{figure}
%
\item
Compute the target Halo orbit in the $\H_1$ family for which the unstable manifold is to
be computed:

The libration point $\Lib_1$ has two pairs of purely imaginary eigenvalues and a
pair of real eigenvalues. The two pairs of purely imaginary eigenvalues correspond
to the 
planar Lyapunov family $\L_1$ and Vertical family $\V_1$
that bifurcate from $\Lib_1$.
Compute
the family $\L_1$, using a standard starting procedure \citep{AUTO81} at the libration
point $\Lib_1$.
The free problem parameters in this continuation are the period $T$ and the unfolding
parameter $\sigma$. Along $\L_1$ two branch points are located; see Fig.~\ref{fig:l1v1}.
Branch switching at the first of these branch points gives the Halo family $\H_1$.
\item
Determine the eigenfunction of a selected Halo orbit:

Select an appropriate periodic orbit in the $\H_1$ family that has one real Floquet
multiplier with
absolute value greater than 1, so that its unstable manifold is two-dimensional.
Compute the corresponding unstable eigenfunction as follows:
Couple the boundary value equations for $\v$ in Eq.~\eqref{eq:eigenfunction}
to those for $\u$ in Eqs.~\eqref{eq:ODE}, \eqref{eq:bc}, \eqref{eq:phase}.
Supplement this extended system by an appropriate continuation equation of the
form of Eq.~\eqref{eq:pseudo}, that also includes $\v$, $\lambda$, and $\rho$,
with $\v$ and $\rho$ initialized to zero, and $\lambda$ initialized to the
desired Floquet exponent as obtained from the decomposition of the Jacobian 
of the collocation equations.  Written out, this gives
\begin{equation}
\begin{split} 
&\dot{\u}(t) = T \f(\u(t),\sigma)~,\\
&\u(1) = \u(0)~,\\
&\int_0^1 \langle \u(t),\dot{\u}_0(t) \rangle dt ~=~0~,\\
&\dot\v(t) = T \f_{\u} \bigl(\u(t),0\bigr) \v(t)+\lambda \v(t)~,\\
&\v(1) = \v(0)~, \\
&\langle \v(0),\v(0) \rangle =\rho~,\\
&\int_0^1
       \langle \w(t)-\w_0 (t), \w_0' (t)\rangle dt
       + \langle \p-\p_0, \p_0' \rangle
       = \Delta s~, \quad \w(t)=(\u(t),\v(t)),\quad \p=(\sigma,\lambda,\rho).\\
\end{split}
\label{eq:eigenfunctionextended}
\end{equation}
In our continuation context we observe that the target period orbit $\u(t)$,
together with $\v(t)\equiv0$ and $\p=(\sigma,\lambda,\rho)=(0,\lambda,0))$,
corresponds to a {\it branch point}, from which a solution family
$(\u,\v;\sigma,\lambda,\rho)$ bifurcates. Along this bifurcating family
the orbit $\u(t)$ remains equal to the target periodic orbit, the Floquet
exponent $\lambda$ and the period $T$ of $\u(t)$ remain constant, while
the unfolding parmameter $\sigma$ remains zero; all to numerical precision.
On the other hand, the Floquet eigenfunction $\v(t)$ becomes nonzero, 
and consequently $\rho$ also becomes nonzero. 
In fact, one continuation step along the bifurcating family would be enough
to obtain the nonzero eigenfunction $\v(t)$. However, for numerical purposes we
typically do a few continuation steps along the bifurcating family until $\rho$ 
is equal to $1$, {\it i.e.}, until the norm of the Floquet eigenfunction equals $1$.
This simple procedure to determine the Floquet eigenfunction fits well into the 
continuation and branch-switching algorithms in AUTO, it is numerically stable, 
and has the additional advantage that the Floquet eigenfunction, once computed, 
can be continued with fixed nonzero $\rho$ and varying energy. In that case 
$\p=(\sigma,\lambda,T)$ in the last constraint of 
Eq.~\eqref{eq:eigenfunctionextended}. This allows the
determination of the eigenfunction $\v$ in highly sensitive cases, for example,
for very large or very small Floquet multipliers.
\end{enumerate}
Once the periodic orbit and its appropriate eigenfunction have been computed,
we can proceed to compute the manifold.
\section{Basic Computation of Unstable Manifolds}
\label{sec:manifold-algorithm}
Here we describe a method for computing unstable manifolds based on
the continuation of orbits that lie in it. These orbits start (as a function of
time) in the linear approximation of the unstable manifold that was computed
above, and end in a section where one of the coordinates is fixed.
Other possibilities to constrain the end point include fixing the integration
time or the arclength; we found the fixed end section to be the most appropriate
here.  All these approaches are robust against sensitive dependence on initial
conditions; see also \citet{krauskopf:07}.

\begin{figure}[htb]
\begin{center}
\includegraphics{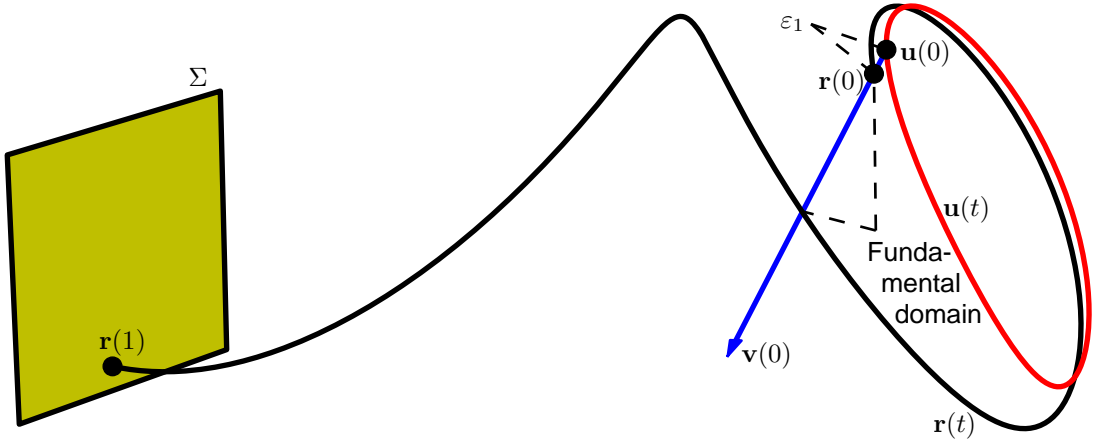}
\end{center}
\caption{Plot of a periodic orbit $\u(t)$ having
Floquet eigenfunction $\v(t)$.
Also shown is the orbit $\r(t)$ in the unstable manifold that starts at the
point $\r(0) = \u(0)+ \eps_1 \v(0)$ and ends in a section $\Sigma$.
Here $\u(t)$ is an $\H_1$ Halo orbit, where
$T=2.3200$, $E=-1.5052$, $\lambda=1.4534$, $\eps_1=0.05$, and
the section $\Sigma$ is at $x=0$. For accuracy the absolute value of $\eps_1$
should be smaller, but for clarity of the diagram it is taken larger here.
}
\label{basicmanifold}
\end{figure}

We now explain the procedure for computing the unstable manifolds of a
given periodic orbit $\u(t)$ in some detail.
We assume that $\u(t)$  has one real Floquet multiplier with
absolute value greater than $1$, with associated eigenfunction $\v(t)$,
computed as described in Sect.~\ref{sec:periodic-orbits}.
The starting data needed are a point on the periodic orbit, namely $\u(0)$,
and the corresponding value of the Floquet eigenfunction, namely $\v(0)$.
For a given appropriately small value of $\eps_1$, take the point
$\r(0) = \u(0) + \eps_1 \v(0)$ as the starting point of an orbit $\r(t)$ in
the unstable manifold; see Fig.~\ref{basicmanifold}.
Here $t$ denotes the time along the orbit $\r(t)$.
Similar to the case of the periodic orbit $\u(t)$ and its eigenfunction $\v(t)$,
where for numerical reasons the time interval was rescaled from $[0,T]$ to $[0,1]$, 
for the non-periodic orbit $\r(t)$ we also rescale time to the unit interval.
Select a section $\Sigma$, for example, at
$x_\Sigma=0$ or $x_\Sigma=-0.25$,
where the orbit $\r(t)$ is to terminate; {\it i.e.}, $\r(1) \in \Sigma$.
(We may allow the orbit $\r(t)$ to cross $\Sigma$ several times, before it actually
terminates in $\Sigma$.)
The part of the manifold to be approximated is then given by the set of orbits
\begin{equation}
\{\r(t)~|~\r(0)=\u(0)+\eps \v(0) ~\text{ and }~ \r(1) \in \Sigma,
  ~\text{ for }~ \eps_1 \le \eps < \eps_2 \}~.
\end{equation}
The range of values of $\eps$, namely, $[\eps_1,\eps_2)$, should be chosen
to correspond to a \emph{fundamental domain}, see Fig~\ref{basicmanifold}.
This ensures that the full manifold
is swept out, at least locally near the periodic orbit, as $\eps$ is allowed to
vary from $\eps_1$ to $\eps_2$.
If the interval $[\eps_1,\eps_2)$ is chosen too small then the manifold would
be incomplete. Taking this interval too large would lead to duplication of orbits,
albeit having different lengths.
A fundamental domain is such that the orbit that starts at $\u(0)+\eps_1 \v(0)$
closely passes the line given by $\u(0)+\eps \v(0)$ again, for the first
time, at $\u(0)+\eps_2 \v(0)$.
For a given value of $\eps_1$, using fundamental properties of the
eigenfunction $\v(t)$ we can compute the corresponding linear approximation of
$\eps_2$ using $\eps_2 = e^{\lambda} \eps_1$.

We also note that if $\eps_1$ is too small in absolute value then the
orbit $\r(t)$ remains close to the periodic orbit for a long time
before it escapes to ultimately reach the section $\Sigma$.
If, on the other hand, $\eps_1$ is too large in absolute value then the
linear approximation of the manifold is no longer accurate.

The boundary value problem for $\r(t)$ is then given by
\begin{equation}
\begin{split}
&\dot\r(t)  ~=~  T_r~\f(\r(t),0)~,\\
&\r(0)   ~=~ \u(0)+\eps \v(0)~,\\
&\r(1)_x ~=~ x_{\Sigma}~,\\
&\int_0^1
       \langle \r(t)-\r_0 (t), \r_0' (t)\rangle dt
       + (\eps-\eps_0) \eps_0' + (T_r-{T_r}_0) {T_r}_0'
       = \Delta s~,\\
\end{split}
\label{eq:sweep}
\end{equation}
where $\sigma=0$ in $\f(\r(t),\sigma)$, and $\Sigma$ denotes the section
$x=x_{\Sigma}$.
As done earlier for the periodic orbit and for its eigenfunction,
the actual integration time $T_r$ of the orbit $\r(t)$
appears explicitly in the differential equation, due to the
above-mentioned scaling of the time $t$.
The pseudo-arclength constraint, the last equation in Eq.~\eqref{eq:sweep}
plays a crucial role in the algorithm. For a given step size $\Delta s$ the
pseudo-arclength constraint ensures in particular that the next computed 
orbit $\r(t)$ is close to $\r_0(t)$ over the whole trajectory, and that 
also $T_r$ is close to ${T_r}_0$.
The step size $\Delta s$ could in principle be constant, but for efficiency
and robustness it is adjusted after each successful continuation step,
depending on the speed of convergence of the Newton/Chord iterations.

Given the point $\u(0)$ on the periodic orbit, and the point
$\r(0) = \u(0)+\eps_1 \v(0)$ in the direction of the unstable manifold
at time $t=0$, the complete procedure to compute the manifold is then
as follows:
\begin{enumerate}
\item
Compute a {\em starting orbit} in the manifold, using continuation to do
the ``time integration''.
More precisely, we use numerical continuation with $T_r$ as free parameter
to compute a ``family'' of solutions, {\it i.e.}, the same trajectory, but
for a set of increasing integration times $T_r$.
The equations used are
\begin{equation}
\begin{split}
\dot\r(t) ~=~  T_r~\f(\r(t),0)~,\\
\r(0)  ~=~ \u(0)+\eps_1 \v(0)~,
\end{split}
\label{eq:grow}
\end{equation}
which correspond to Eq.~\eqref{eq:sweep}, but without the end point constraint.
Note that $\eps$ is fixed at $\eps=\eps_1$ in each continuation step of this
starting procedure.
The continuation is stopped when $\r(1)$ intersects the plane $\Sigma$.
(This termination point need not necessarily be the first such intersection.)
The starting orbit in this continuation is the constant solution
$r(t)\equiv \u(0)+\eps_1 \v(0)$, with $T_r=0$.
Although it may appear inefficient to do a time-integration by continuation,
this approach has the advantage that it fits very well into the continuation
framework of the algorithms in this paper.
\item
Given an orbit that ends in $\Sigma$, as computed above, the unstable
manifold is then approximated by further continuation of this orbit,
now using the boundary value problem in Eq.~\eqref{eq:sweep},
with the end point $\r(1)$ constrained to remain in $\Sigma$,
and with $\eps$ and the integration time $T_r$ allowed to vary.
If $\eps$ varies along a full fundamental domain,
the continuation would sweep out the full unstable manifold,
limited only by the termination condition.
However, the pseudo-arclength constraint in Eq.~\eqref{eq:sweep} limits
the size of the change in the orbit \emph{and} in the the parameters
in any continuation step.
Hence the full unstable manifold is only obtained if no
``obstacles'' are encountered, where for example $T_r$ goes to
infinity. As will be seen in the next section, these obstacles can
in fact be of much interest.
\end{enumerate}
\section{Example Computations of Unstable Manifolds and Connecting Orbits to Invariant Tori}
\label{sec:manifold-examples}

\begin{figure}[htb]
\begin{center}
\begin{picture}(460,280)
\put(0,0){\includegraphics[scale=0.33,clip,trim=0mm 20mm 0mm 0mm]{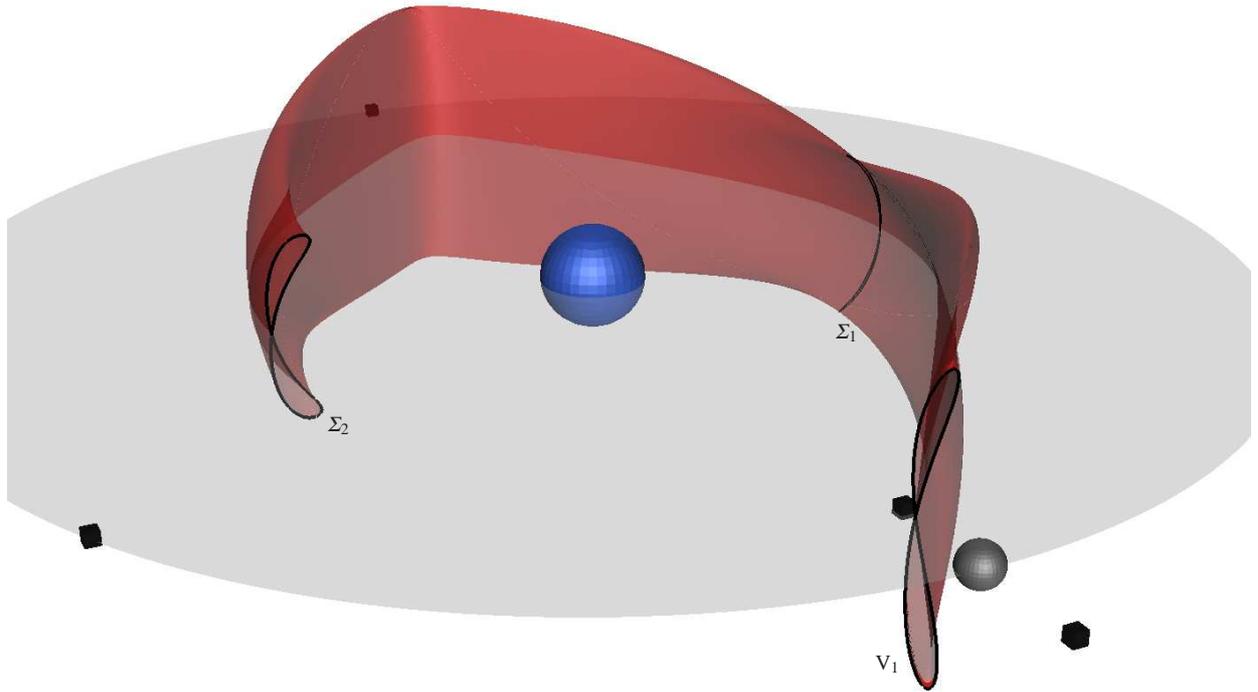}}
\put(325,14){$\V_1$}
\put(120,104){$\Sigma_2$}
\put(310,139){$\Sigma_1$}
\end{picture}
\end{center}
\caption{The unstable manifold of a periodic orbit from the $\V_1$ family.
The periodic orbit, labeled $\V_1$, is located on the $\Lib_1$ side of the Moon,
 with period $T=3.7700$ and energy $E=-1.5164$.
The terminating plane $\Sigma$ is located at $x_{\Sigma}=0$. The first intersection 
of the manifold with $\Sigma$ is indicated by the curve labeled $\Sigma_1$, and the 
second intersection, where the manifold computation is terminated, is labeled $\Sigma_2$.}
\label{fig:manV1}
\end{figure}

\begin{figure}[htb]
\begin{center}
\begin{picture}(460,230)
\put(0,0){\includegraphics[scale=0.33,clip,trim=0mm 40mm 0mm 20mm]{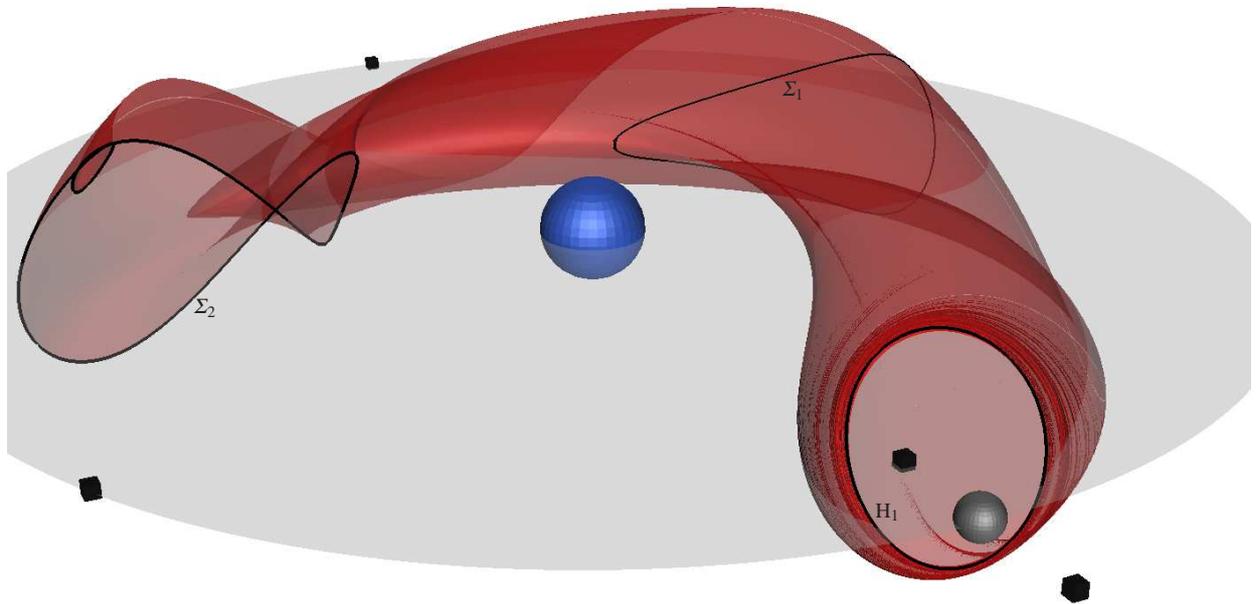}}
\put(325,35){$\H_1$}
\put(70,112){$\Sigma_2$}
\put(290,194){$\Sigma_1$}
\end{picture}
\end{center}
\caption{
The unstable manifold of a periodic orbit (labeled $\H_1$) from the $\H_1$ family.
The periodic orbit, which is on the $\Lib_1$ side of the Moon,
has period $T=2.5152$ and energy $E=-1.5085$.
The terminating plane $\Sigma$ is located at $x_{\Sigma}=-0.25$, and the first and second intersections
of the manifold with $\Sigma$ are labeled $\Sigma_1$ and $\Sigma_2$.}
\label{fig:manH1}
\end{figure}

%
\begin{figure}[htb]
\begin{center}
\includegraphics[scale=0.33,clip,trim=0mm 25mm 0mm 4mm]{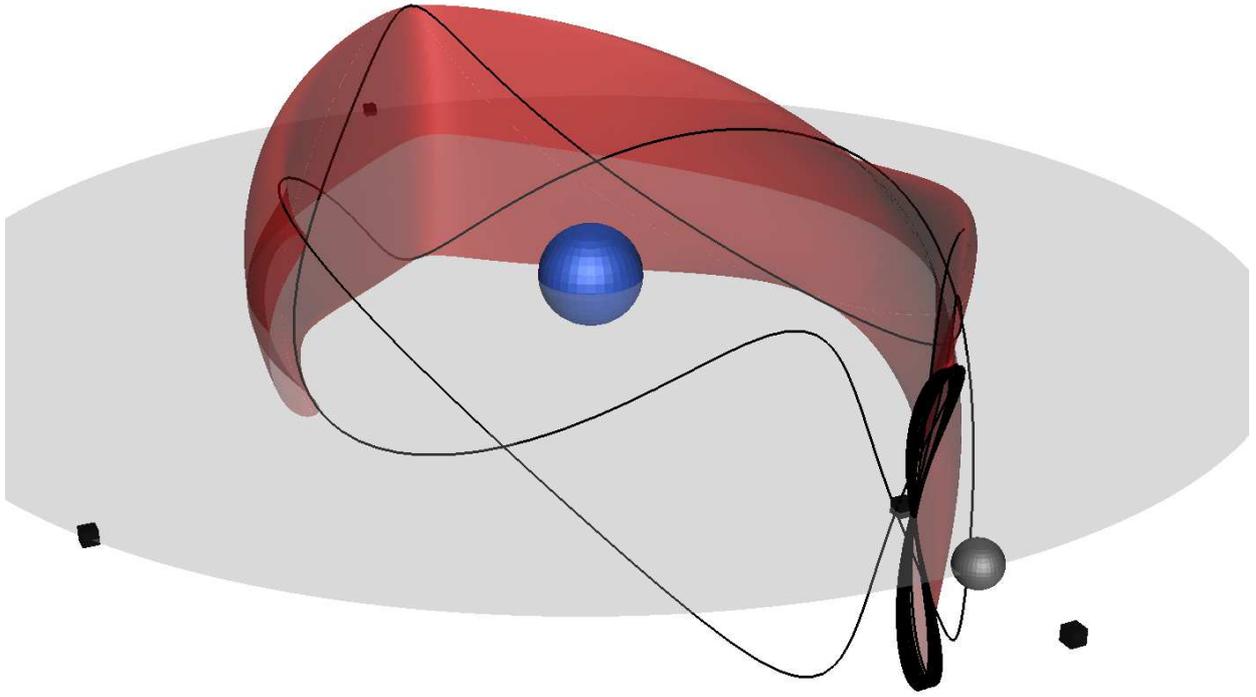}
\end{center}
\caption{
Part of the unstable manifold of the $\V_1$ periodic orbit from Fig.~\ref{fig:manV1},
together with a superimposed longer orbit in this manifold, computed as
described in the text.  }
\label{fig:man+orbit-V1}
\end{figure}
%
\begin{figure}[htb]
\begin{center}
\includegraphics[scale=0.33,clip,trim=0mm 25mm 0mm 4mm]{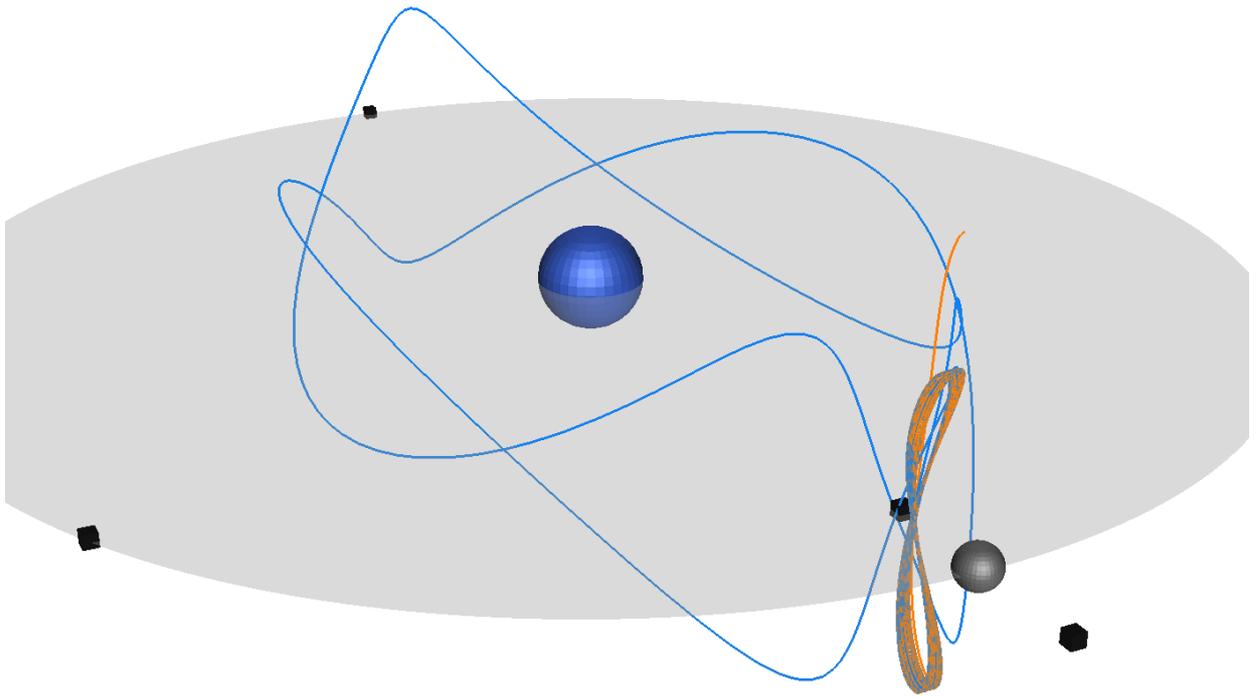}
\end{center}
\caption{
A separate view of the longer orbit from Fig.~\ref{fig:man+orbit-V1},
which winds around a torus near the original $\V_1$ periodic orbit.
The orbit ultimately returns to the plane $\Sigma$
that is located at $x_{\Sigma}=0.5$.
We remark that in this and subsequent figures
the orbit color changes from a deep sky blue (RGB code
(0,0.5,1)) via gray (code (0.5,0.5,0.5)) to dark orange (code
(1.0,0.5,0)) as the scaled time increases from 0 to 1.
}
\label{fig:orbit-V1}
\end{figure}

\begin{figure}[htb]
\begin{center}
\includegraphics[scale=0.33,clip,trim=5mm 20mm 0mm 3mm]{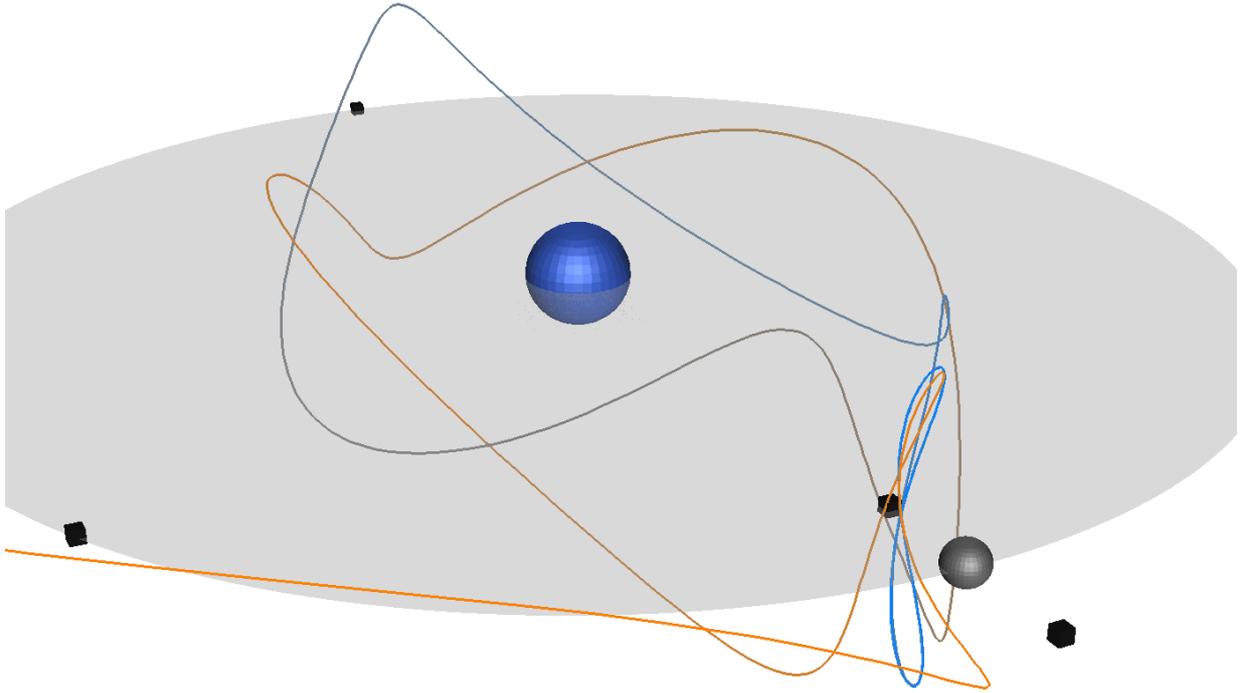}
\end{center}
\caption{An orbit obtained by initial value integration for the same
value of $\eps$, within numerical precision, as used for the orbit
in Fig.~\ref{fig:orbit-V1}. Note that the orbit approaches but then
diverts from the torus seen in Fig.~\ref{fig:orbit-V1}, without
winding around it.}
\label{fig:shooting}
\end{figure}

\begin{figure}[htb]
\begin{center}
\includegraphics[scale=0.33,clip,trim=5mm 27mm 0mm 3mm]{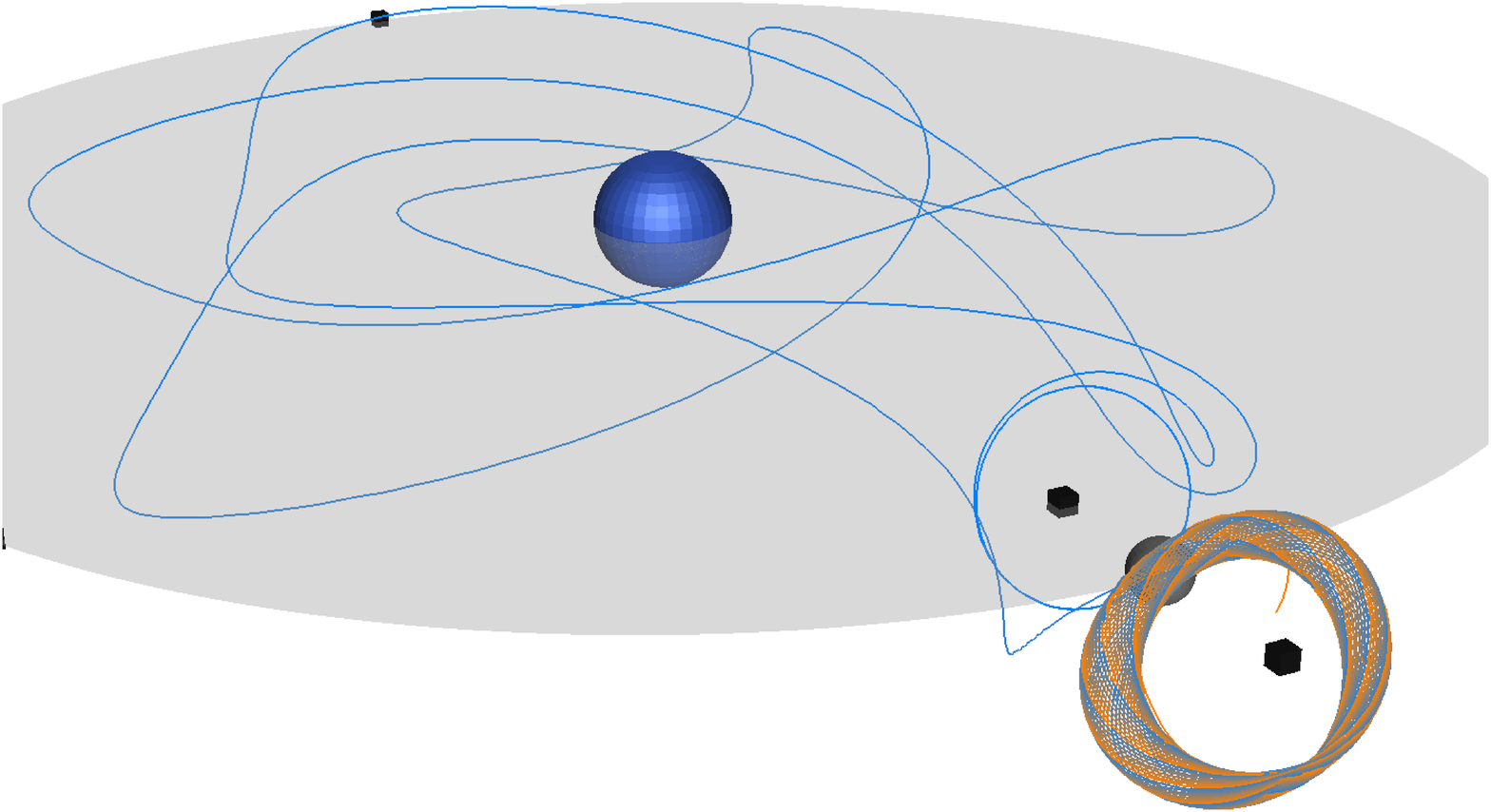}
\end{center}
\caption{
Continuation of a longer orbit within the unstable manifold of
a Halo orbit with $E=-1.5532$
with the plane $\Sigma$ located at $x_{\Sigma}=1.02$, that is, on the
$\Lib_2$-side of the Moon. This orbit winds around a torus near an orbit
from the Halo family $\H_2$, before returning to and ending in the plane
$\Sigma$.
}
\label{fig:H1-to-H2-torus}
\end{figure}
\afterpage{\clearpage}

In this section we illustrate the basic boundary value technique described
above by computing unstable manifolds of the Vertical family $\V_1$ and the
Halo family $\H_1$. For certain ranges of their energy (the thin
solid curves in Fig.~\ref{fig:families}) and period these
families have periodic orbits with one real, positive Floquet multiplier
outside the unit circle, so that their unstable manifold is indeed
two-dimensional.
The algorithm applies in principle also to higher-dimensional unstable manifolds,
by using a selected unstable Floquet multiplier, typically the largest one,
and its associated eigenfunction. As already mentioned, the algorithm also
applies to stable manifolds. However, here we restrict attention to examples
with two-dimensional unstable manifolds.

Fig.~\ref{fig:manV1} shows the computed unstable manifold of an orbit from the
Vertical family $\V_1$. Here the section $\Sigma$ is taken at $x_\Sigma=0$.
As seen in the figure, the orbits in the manifold terminate in $\Sigma$
at their second intersection with this plane.
Note also that the intersection curve of the manifold with the section
$\Sigma$ has a similar shape as the figure-eight Vertical orbit from which
the manifold originates.

Fig.~\ref{fig:manH1} shows the computed unstable manifold of a Halo orbit
from the $\H_1$ family. The plane $\Sigma$ is located at $x_\Sigma=-0.25$.
Note that the manifold changes shape as it propagates.
Here, as well as for the unstable manifold of the $\V_1$ orbit,
there is no contradiction in the fact that the cross section of the manifold
with the plane $\Sigma$ is a self-intersecting curve, since we are viewing
a projection of the manifold from $\mathbb{R}^6$ (spatial and velocity variables)
into $\mathbb{R}^3$ (spatial variables only).

If the section is taken at certain other value-ranges of $x_{\Sigma}$ then
one may encounter obstacles that prevent the continuation to cover the full
manifold. More specifically, the initial value of the orbits, as determined
by the value of $\eps$, does not cover the entire fundamental domain, 
but approaches a particular value where $T_r\to\infty$.
One possible scenario is that the orbit picks up additional loops around 
a torus-like object near its end point, if such an object exists, before 
returning to and ending in the plane $\Sigma$.
This computational phenomenon is not exceptional in the CR3BP; in fact it is
easy to find specific examples. One such instance is given in
Fig.~\ref{fig:man+orbit-V1},
which again shows the unstable manifold of the Vertical orbit in
Fig.~\ref{fig:manV1}, but with a superimposed longer orbit that lies in
the same unstable manifold.  This longer orbit alone is shown separately in
Fig.~\ref{fig:orbit-V1}. As is evident from Fig.~\ref{fig:orbit-V1}, the
orbit returns to a neighborhood of the original Vertical orbit, where it winds
around a torus-like object (a \emph{quasi-Vertical} orbit).

More specifically, the phenomenon shown in Fig.~\ref{fig:man+orbit-V1} and
Fig.~\ref{fig:orbit-V1} results from ``growing'' (see step 1 at the
end of Sect.~\ref{sec:manifold-algorithm})
a longer initial orbit in the
unstable manifold of the $\V_1$ periodic orbit, with $\eps$ fixed,
until it intersects a plane $\Sigma$ located at $x_{\Sigma}=0.5$.
Continuing this orbit
with the end point constrained to remain in $\Sigma$, and with $\eps$
and $T_r$ allowed
to vary, then results in it winding around a torus near the original $\V_1$
periodic orbit, before returning to and terminating in the plane $\Sigma$.
In contrast, trying to obtain such an orbit directly by shooting techniques 
appears to be difficult, as such orbits approach but then divert from the 
torus; see Fig.~\ref{fig:shooting}.
Indeed, for the orbits to agree it may be necessary to compute the BVP
orbits to higher precision, {\it i.e.}, more than the double precision
used in AUTO. In particular, this would yield the initial value to
higher precision.
Likewise, the initial value integration, starting from such a more accurate
initial value, may also need to be in high precision arithmetic.

A similar result is shown in Fig.~\ref{fig:H1-to-H2-torus}
for an extended orbit from the $\H_1$ manifold shown in Fig.~\ref{fig:manH1}.
Here an orbit in the unstable manifold of the $\H_1$ periodic orbit is grown
until it reaches a plane $\Sigma$ located at $x_{\Sigma}=1.02$, {\it i.e.}, on
the $\Lib_2$-side of the Moon. Constraining the end point of this longer orbit
to remain in $\Sigma$, and allowing $\eps$ to vary, then
results in the end portion of the orbit to wind around a torus
near an orbit from the Halo family $\H_2$ (a \emph{quasi-Halo} orbit),
as shown in Fig.~\ref{fig:H1-to-H2-torus}. Meanwhile the initial value $\r(0)$
of the orbit approaches a point in the fundamental domain, without fully
covering that domain.
Ultimately the orbit returns to the plane $x_{\Sigma}=1.02$, as required by the
computational set-up.

Similarly, in Fig.~\ref{fig:H1-to-H1-torus}
an extended orbit from the $\H_1$ manifold is grown
until it reaches a plane $\Sigma$ located at $x_{\Sigma}=0.6$, after going
around the Earth once. Constraining the end point of this orbit
to $\Sigma$, while $\eps$ and the integration time $T_r$ are allowed to vary, results
in an orbit that winds around a quasi-Halo orbit near $\H_1$.

As a final example in this section we compute an initial orbit in the
unstable manifold of an $\H_1$ Halo orbit, but now in the part of the
unstable manifold on the side of the Moon.
The terminating plane
$\Sigma$ was taken at $x_{\Sigma}=1.02$. Further continuation of this
orbit with $x_{\Sigma}$ constrained to remain in $\Sigma$, and with
$\eps$ allowed to vary, results in the end portion of this orbit winding around
a quasi-Halo orbit near an $\H_2$ Halo orbit,
as shown in Fig.~\ref{fig:H1-to-H2-Direct},
before returning to and ending in $\Sigma$.

In this section we have computed the unstable manifolds of some orbits in the
$\H_1$ Halo and $\V_1$ Vertical families. For certain values of $\eps$
we find what appear to be heteroclinic connecting orbits from the original
periodic orbit to an invariant torus. These connecting orbits appear to be 
more difficult to obtain with shooting techniques.
This was seen in Fig.~\ref{fig:shooting}, where shooting failed to reveal 
the invariant torus, when starting near the numerical value of $\eps$ 
for which the heteroclinic connection was found in the BVP approach.
This is in part due to the fact that the tori are evidently unstable, with 
saddle-type stability since the connecting orbit approaches it, but ultimately
also leaves the neighborhood of the torus.

The next objective is to continue the orbit-to-torus connections as the energy
$E$ is allowed to vary, and with it the base periodic orbit itself.
In the following section we explain the computational set-up, while examples
are given in Sect.~\ref{sec:connecting-orbits-examples}.
In particular we will see in Sect.~\ref{sec:connecting-orbits-examples} how
such continuation can lead to homoclinic and heteroclinic connecting orbits
between periodic orbits in the same or different families.

\begin{figure}[htb]
\begin{center}
\includegraphics[scale=0.33,clip,trim=5mm 42mm 0mm 32mm]{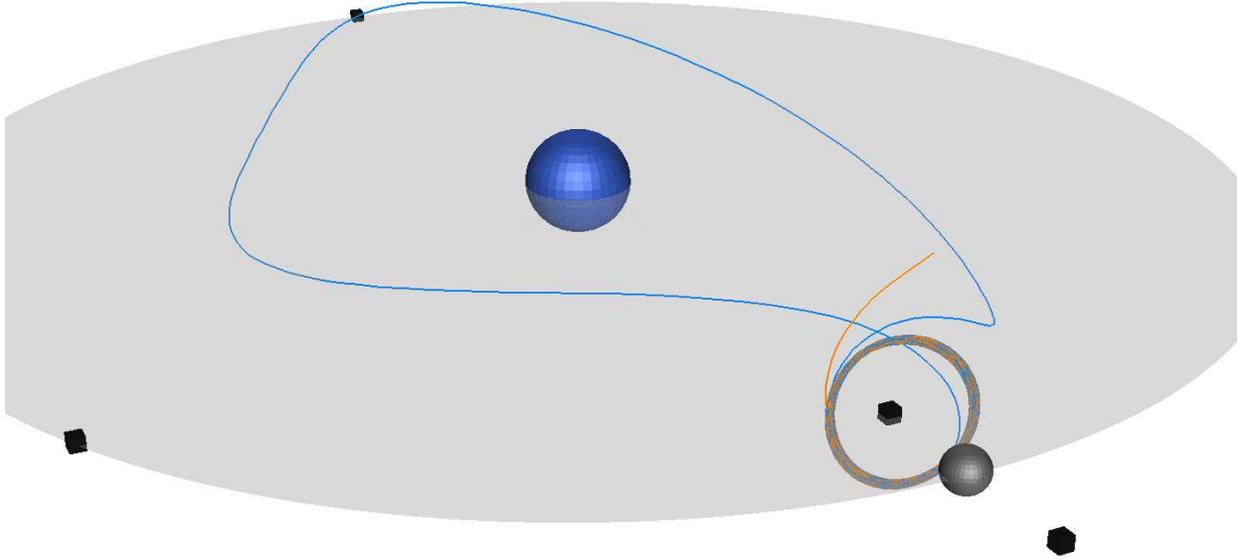}
\end{center}
\caption{
Continuation of a longer orbit within the unstable manifold of
a Halo orbit with $E=-1.5631$
with the plane $\Sigma$ located at $x_{\Sigma}=0.6$.
This orbit winds around a torus near the orbit
from the Halo family $\H_1$ where it started from,
before returning to and ending in the plane $\Sigma$.
}
\label{fig:H1-to-H1-torus}
\end{figure}

\begin{figure}[htb]
\begin{center}
\includegraphics[scale=0.33,clip,trim=0mm 66mm 0mm 50mm]{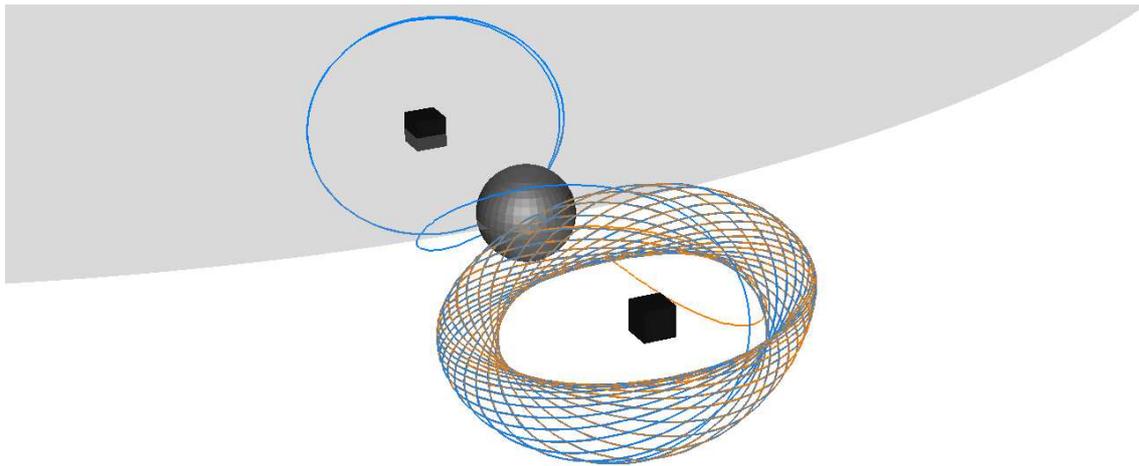}
\end{center}
\caption{
An orbit in the Moon-side part of the unstable manifold of an $\H_1$ periodic
orbit. This orbit is found by fixing its end point to remain in a plane $\Sigma$
located at $x_{\Sigma}=1.02$, that is, on the $\Lib_2$-side of the Moon.
The orbit winds around a torus near an orbit from the Halo family $\H_2$, before
returning to and ending in the plane $\Sigma$.
}
\label{fig:H1-to-H2-Direct}
\end{figure}
\section{Continuation of Connecting Orbits}
\label{sec:connecting-orbits}

The examples in the preceding section provide evidence for the existence
of orbits in the unstable manifold of certain periodic orbits that connect to
toroidal objects in phase space. These connections exist for specific initial
points in the fundamental domain, {\it i.e.}, along the line $\u(0) + \eps \v(0)$,
$\eps \in [\eps_1,\eps_2)$.
The connecting orbits are generated numerically as the initial value of the orbits
approaches a final point inside the fundamental domain, as the integration time $T_r$ 
increases and additional windings are generated near the end of the orbit.
As a consequence the fundamental domain is not fully covered, although this can be 
remedied by repeating the entire sequence of steps in the algorithm, starting from 
a value of $\eps$ that corresponds to a point in the part of the fundamental domain 
that is not covered.

The question arises as to what happens to the periodic-orbit-to-torus connections
if the periodic orbit is varied along one of the primary or secondary families.
One approach would be to repeat the computations
for each one of a sequence of periodic orbits along a given family, {\it e.g.}, along the
Halo family $\H_1$. However, here also, continuation is a more effective tool
that allows interesting connecting orbits to be detected easily
along the continuation path.

We are then led to consider the following approach:
Collect the equations that define the periodic orbit
(Eqs.~\eqref{eq:ODE},~\eqref{eq:bc}, \eqref{eq:phase}),
the equations defining the Floquet multiplier and eigenfunction
(Eq.~\eqref{eq:eigenfunction}),
and the equations for the orbit in the unstable manifold
(Eq.~\eqref{eq:sweep}).
For clarity, this complete set of coupled equations is reproduced in
Eq.~\eqref{eq:All}.
\begin{equation}
\begin{split}
&\dot{\u}(t) ~=~ T~\f(\u(t),\sigma)~, \\
&\u(1) ~=~ \u(0)~, \\
&\int_0^1 \langle \u(t),\dot{\u}_0(t) \rangle dt ~=~0~,\\
&\dot\v(t) = T \f_{\u} \bigl(\u(t),0\bigr) \v(t)+\lambda \v(t)~,\\
&\v(1) = \v(0)~, \\
&\langle \v(0),\v(0) \rangle =\rho~,\\
&\dot\r(t)  ~=~  T_r~\f(\r(t),0)~,\\
&\r(0) ~=~ \u(0)+\eps \v(0)~,\\
&\r(1)_x ~=~ x_{\Sigma}~,\\
&\int_0^1
       \langle \w(t)-\w_0 (t), \w_0' (t)\rangle dt
       + \langle \p-\p_0, \p_0'\rangle
       ~=~ \Delta s~, \quad \w(t)=(\u(t),\v(t),\r(t)),\quad \p=(T,\sigma,\lambda,\eps)
\label{eq:All}
\end{split}
\end{equation}
The last constraint shown in Eq.~\eqref{eq:All} is the suitably expanded version of the
pseudo-arclength constraint (Eq.~\eqref{eq:pseudo}) that defines the
continuation step size.

Recall that each of the vectors $\u$, $\v$, and $\r$, is six-dimensional.
A simple count shows that Eq.~\eqref{eq:All} represents a system
of $18$ ODEs, subject to a total of $21$ boundary and integral constraints,
not counting the pseudo-arclength constraint. Generically, the continuation of
a solution to Eq.~\eqref{eq:All} then requires four free scalar variables.
The appropriate choice of these parameters is $T$, $\sigma$, $\lambda$, and
$\eps$, where $T$ is the period of $\u(t)$, $\sigma$ is the unfolding parameter,
$e^{\lambda}$ is the Floquet multiplier, and $\eps$ corresponds to the step size
in the direction of the unstable manifold.

The BVP in Eq.~\eqref{eq:All} does not directly define a connection from a periodic
orbit to a torus, but a connection from a periodic orbit to a
section $\Sigma$ for a sufficiently large, fixed, value of the integration time $T_r$ of
$\r(t)$. The torus is then indirectly continued. This indirect
approach is very useful as the direct approach is
known to be of considerable algorithmic complexity, as addressed for
example in
\citet{DieciLorenzRussell,DieciLorenz,EdohRussellSun,henderson2002,
schilder05,olikara2010}.

The strategy of indirect continuation
is somewhat analogous to the computation of a simple homoclinic orbit, {\it i.e.},
an orbit in phase space that approaches a given saddle equilibrium in both
positive and negative time. The continuation of such a homoclinic orbit in
two parameters can be directly formulated as a boundary value problem with asymptotic
boundary conditions that compute the saddle point and its relevant
eigenspaces. In fact, this approach has been used very effectively for the
continuation of homoclinic and heteroclinic orbits, including higher co-dimension
orbits, both for orbits homoclinic or heteroclinic to equilibria or to periodic
orbits \citep{CKS:96}.
However, a generic homoclinic orbit can also be very effectively approximated indirectly
by a periodic orbit of high period, which renders the 2-parameter continuation
of a homoclinic orbit as simple as the continuation of a periodic
orbit (using the techniques given in Sect.~\ref{sec:periodic-orbits}) where
the period $T$ is fixed at a sufficiently large value.

In the next section we demonstrate the use of the indirect
periodic-orbit-to-torus approach described above, by applying it to the
connecting orbits from a Halo orbit in the $\H_1$ family to a torus, as
initially computed in Sect.~\ref{sec:manifold-examples}.
We will see that these continuation calculations can lead to approximate
connecting orbits from the base periodic orbit to resonant periodic
orbits and to other periodic orbits. The detection of such
approximate connections can be refined, which we do not do in this paper,
to provide more accurate approximations to such connecting orbits.
These, in turn, can then be continued directly, using known algorithms
for this purpose, as described, for example, in \citet{CKS:96,dkkv1,dkkv2},
and as implemented in AUTO.

\section{Showcasing Three Families of Connecting Orbits}
\label{sec:connecting-orbits-examples}

In this section we describe the results of the further continuation of orbits
that connect a periodic orbit to a torus, as initially found in
Sect.~\ref{sec:manifold-examples}.
This is done using the $18$-dimensional system in Eq.~\eqref{eq:All}, which
was discussed in the preceding section.

\begin{figure}[htbp]
\begin{center}
\subfloat[]{
\includegraphics[scale=0.26,clip,trim=70mm 20mm 40mm 20mm]{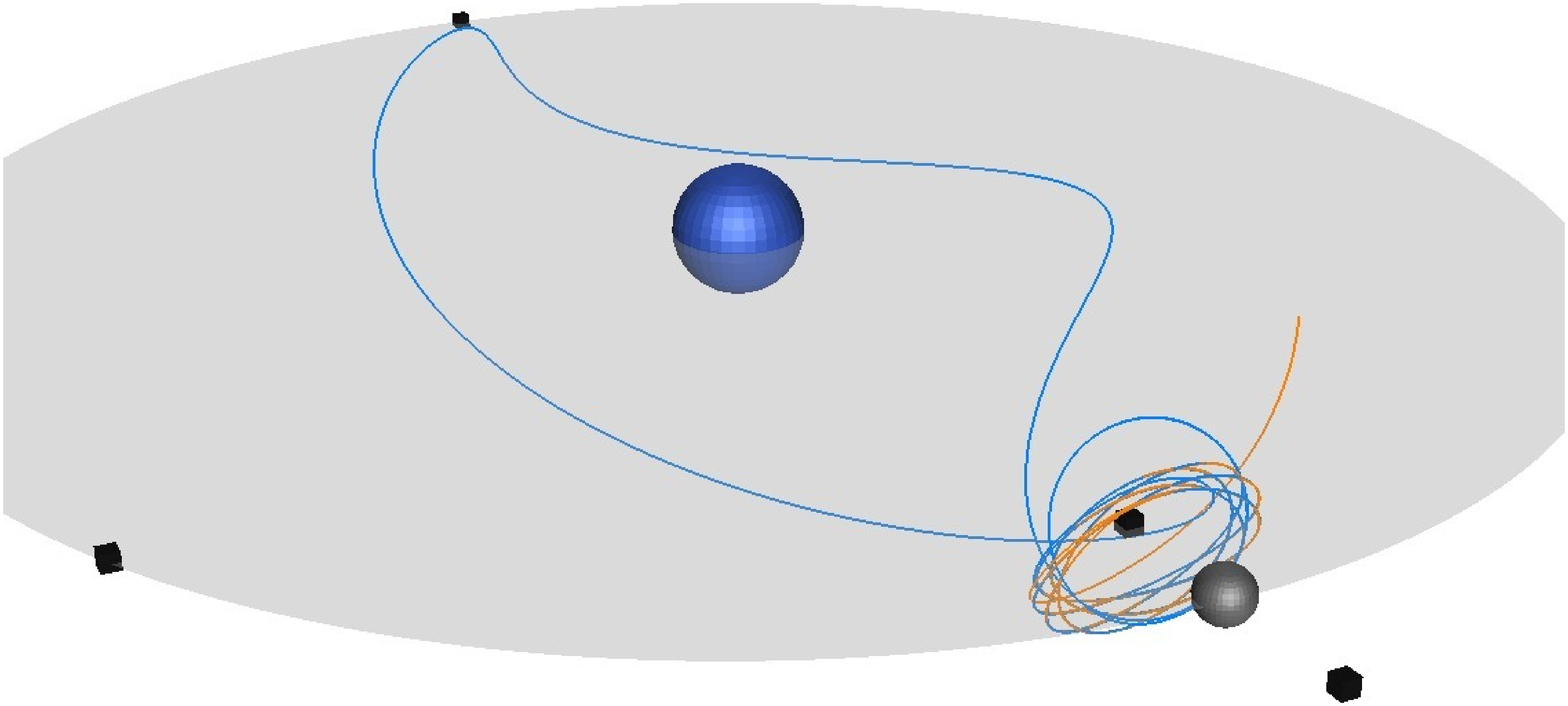}
}
\subfloat[]{
\includegraphics[scale=0.26,clip,trim=70mm 20mm 40mm 20mm]{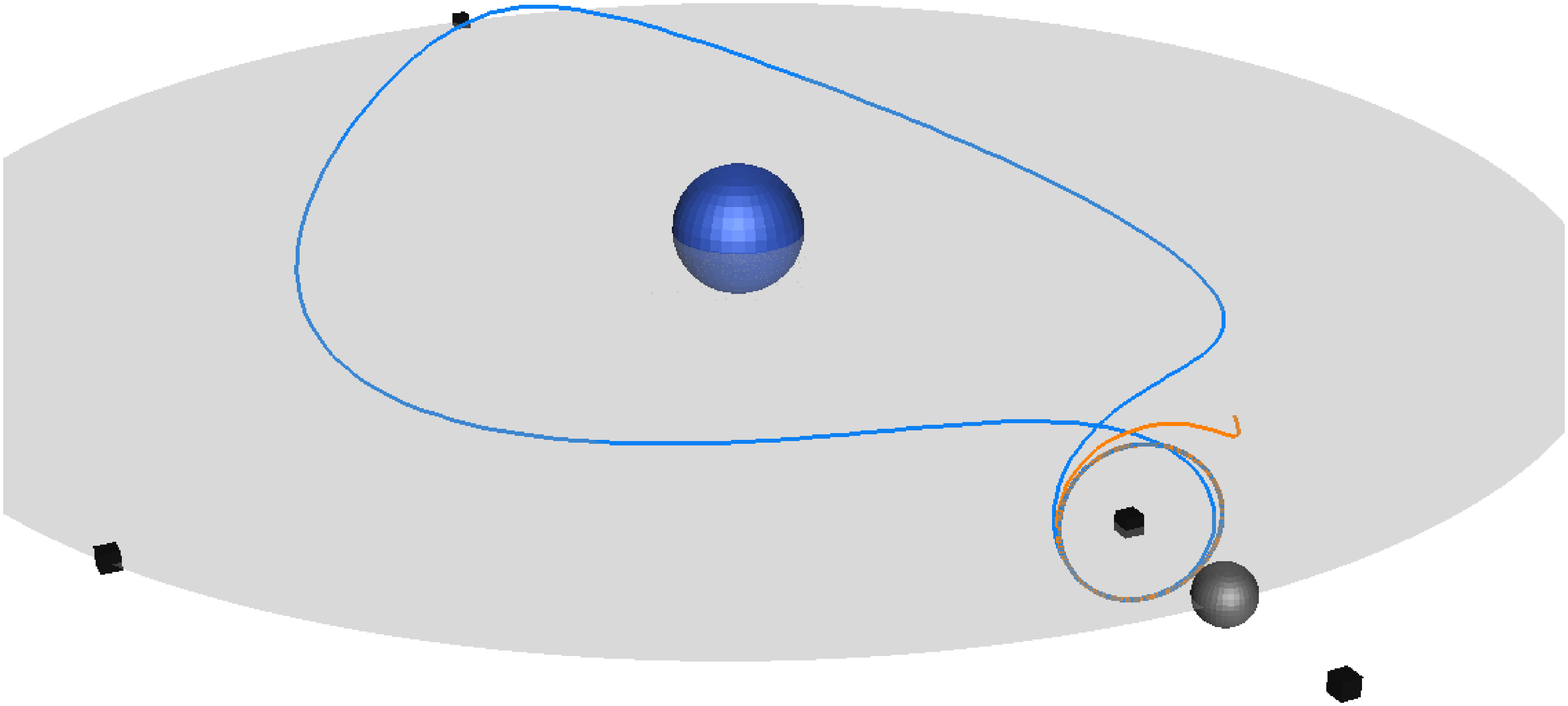}
}

\subfloat[]{
\includegraphics[scale=0.26,clip,trim=70mm 20mm 40mm 20mm]{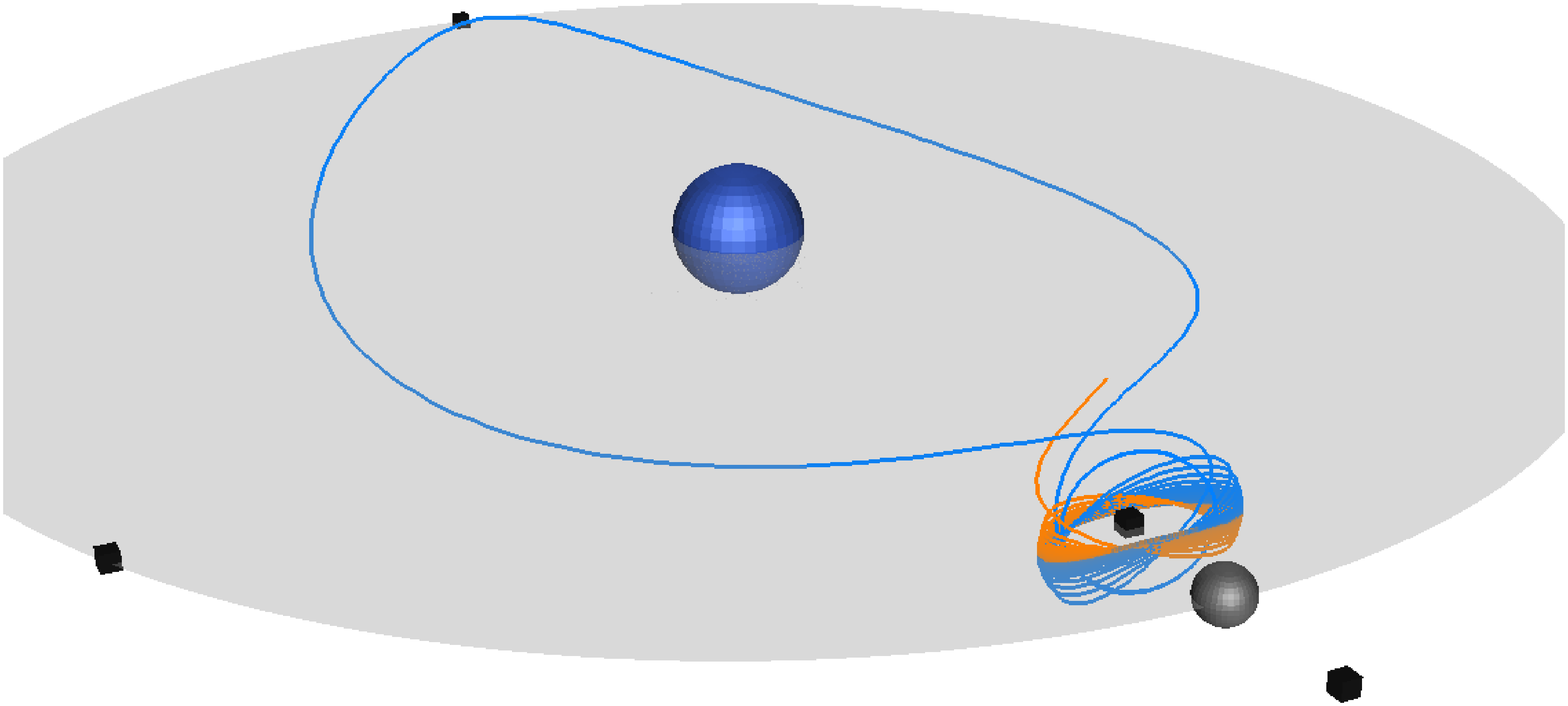}
}
\subfloat[]{
\includegraphics[scale=0.26,clip,trim=70mm 20mm 40mm 20mm]{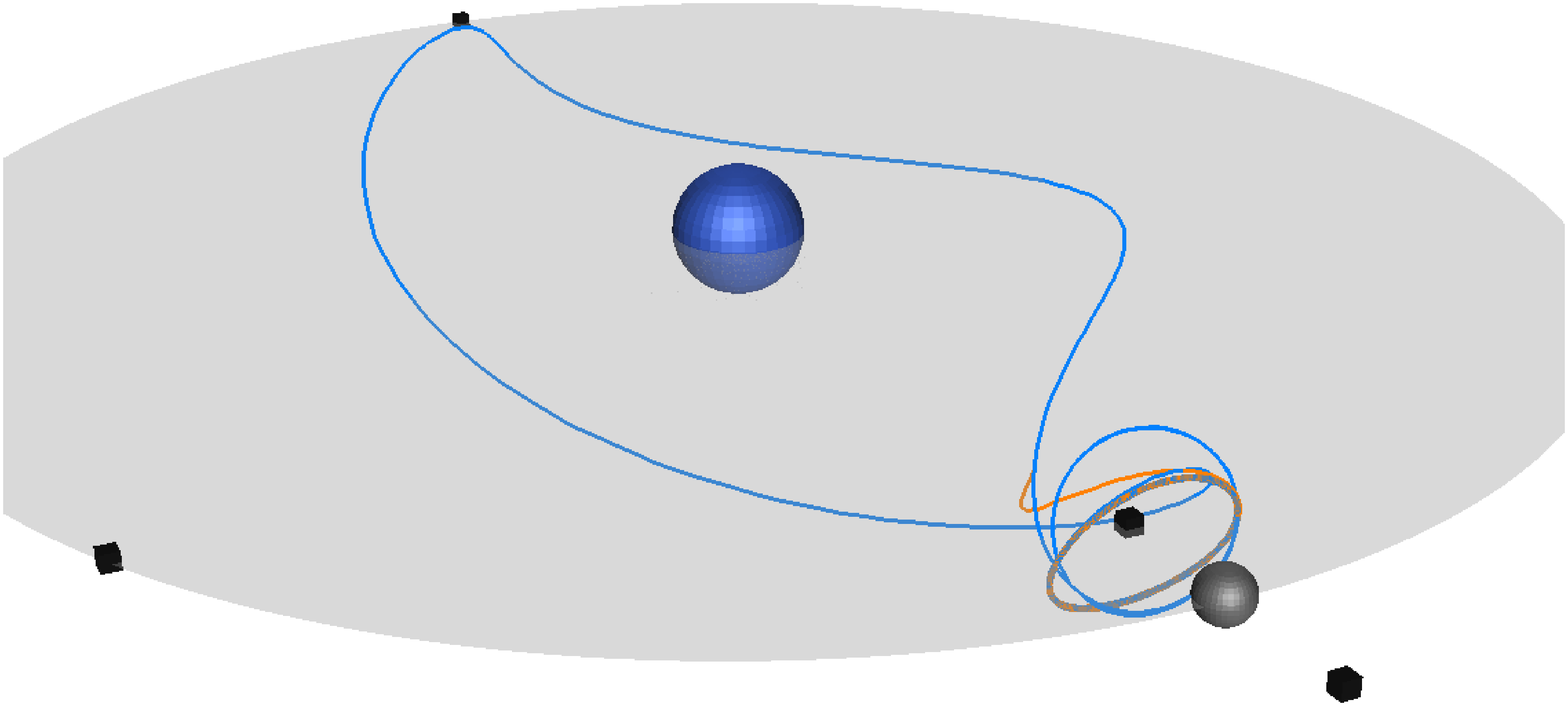}
}

\vskip-0.1truein
\subfloat[]{
\includegraphics{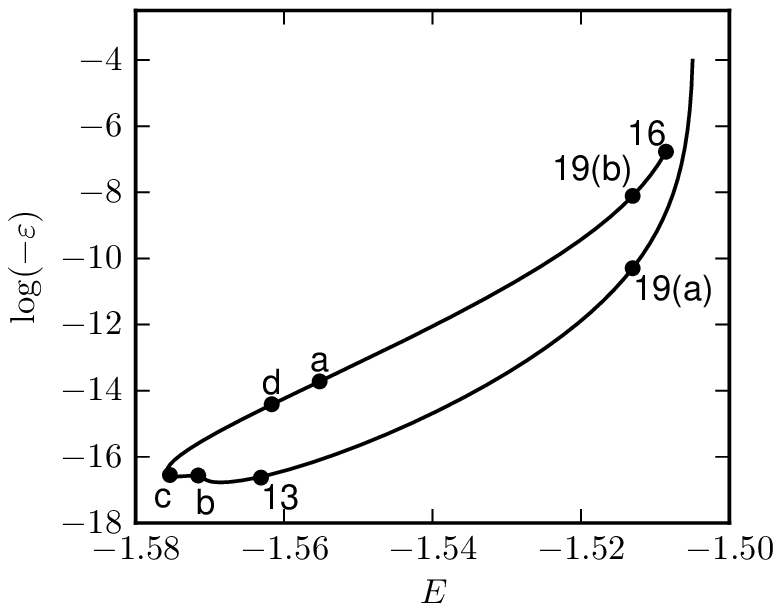}
}
\subfloat[]{
\includegraphics{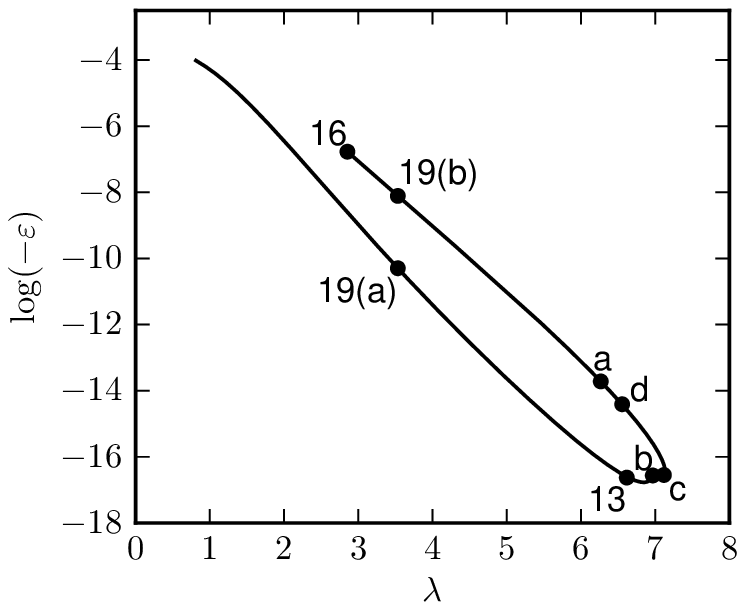}
}
\end{center}
\vskip-0.25truein
\caption{
A family of connecting orbits from $\H_1$ Halo orbits to tori,
making one loop around the Earth.
The terminating plane $\Sigma$ is located at $x_{\Sigma}=0.6$.
(a)--(d): Connecting orbits from
$\H_1$ Halo orbits to (a) a 5:1 resonant orbit near the
corresponding southern Halo orbit, (b) approximately the originating $\H_1$
Halo orbit, (c) an $\L_1$ planar orbit, and (d) a southern $\H_1$ orbit.
(e): Energy of the Halo orbit versus $\log(-\eps)$ along the continuation path.
The continuation terminates in the top right corner
when the Floquet multipliers of the $\H_1$ orbit
become complex, and at label 16 in the heteroclinic
connection shown in Fig.~\ref{fig:H1-to-Axial}.
Labels a--d,~\ref{fig:H1-to-H1-torus},~\ref{fig:H1-to-Axial},
and~\ref{fig:H1-poincare}(a) and (b)
correspond to the values in the panels above and in
Figs.~\ref{fig:H1-to-H1-torus},~\ref{fig:H1-to-Axial}, and~\ref{fig:H1-poincare}.
(f): The logarithm of the unstable Floquet multiplier of the Halo orbit versus $\log(-\eps)$ along
the continuation path.
}
\label{fig:H1-1Loop}
\end{figure}


As a first example we start from the orbit shown in Fig.~\ref{fig:H1-to-H1-torus}
which connects the $\H_1$ Halo orbit with period $T=2.5152$ and energy $E=-1.5164$
to a quasi-Halo toroidal orbit near the originating Halo orbit,
after making one loop around
the Earth. The continuation procedure using the $18$-dimensional ODE in Eq.~\eqref{eq:All}
allows the energy to change, and with it the $\H_1$ orbit
itself, the connecting orbit, and the torus it approaches.
Interesting transitions are encountered along the continuation path, namely,
in the way the changing connecting orbit winds around the changing torus.

Examples are shown in Fig.~\ref{fig:H1-1Loop}. Panel
(a) shows the connecting orbit evidently approaching a 5:1 resonant
orbit. In panel (b) the quasi-Halo orbit has shrunk so it
can not be visually distinguished from the periodic Halo orbit in whose unstable manifold it lies, that is,
the orbit appears to be homoclinic to the $\H_1$ orbit. Such an orbit could be interesting in
space mission design, allowing for occasional large spatial excursions from an $\H_1$ orbit
for negligible energy cost. Panel (c) shows
a heteroclinic connection between an $\H_1$ orbit and the planar
$\L_1$ orbit with the same energy $E=-1.5754$. Recall that the $\L_1$ family
is the Lyapunov family that bifurcates from the libration point
$\Lib_1$, as seen in Fig.~\ref{fig:families}. This heteroclinic connection
could be used for a low cost transfer orbit between the $\H_1$ and
$\L_1$ orbits with energy $E=-1.5754$, even though these orbits are different
in their dynamical properties (see Figs.~\ref{fig:l1v1} and~\ref{fig:h1a1})
and far apart in the bifurcation diagram of the orbits (see Fig.~\ref{fig:families}).
Panel (d) is similar to panel (b), but now the Halo
orbit is connected to its corresponding southern Halo orbit, which is the mirror
image in the $z=0$ of the $\H_1$ orbit. That is, the connection appears heteroclinic
instead of homoclinic.
The energies of these connecting orbits and the periods of the periodic orbits
that they connect to, for these and other connecting orbits in this paper,
are given in Table~\ref{tab:Erg-Typ-Per-Flq}.

\begin{figure}[Htbp]
\begin{center}
\subfloat[]{
\includegraphics[scale=0.33,clip,trim=0mm 10mm 0mm 10mm]{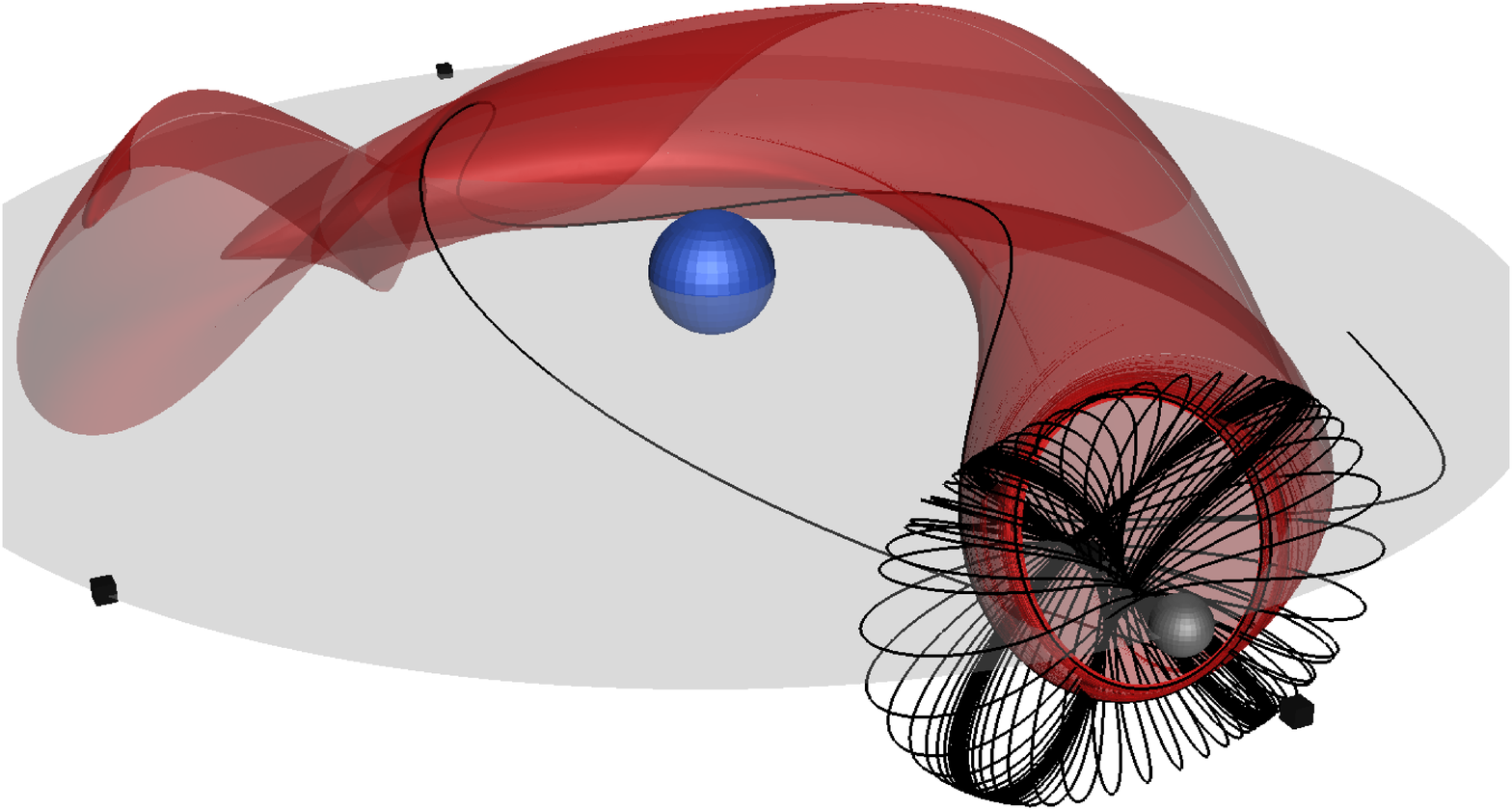}
}
\end{center}
\begin{center}
\subfloat[]{
\includegraphics[scale=0.33,clip,trim=0mm 10mm 0mm 30mm]{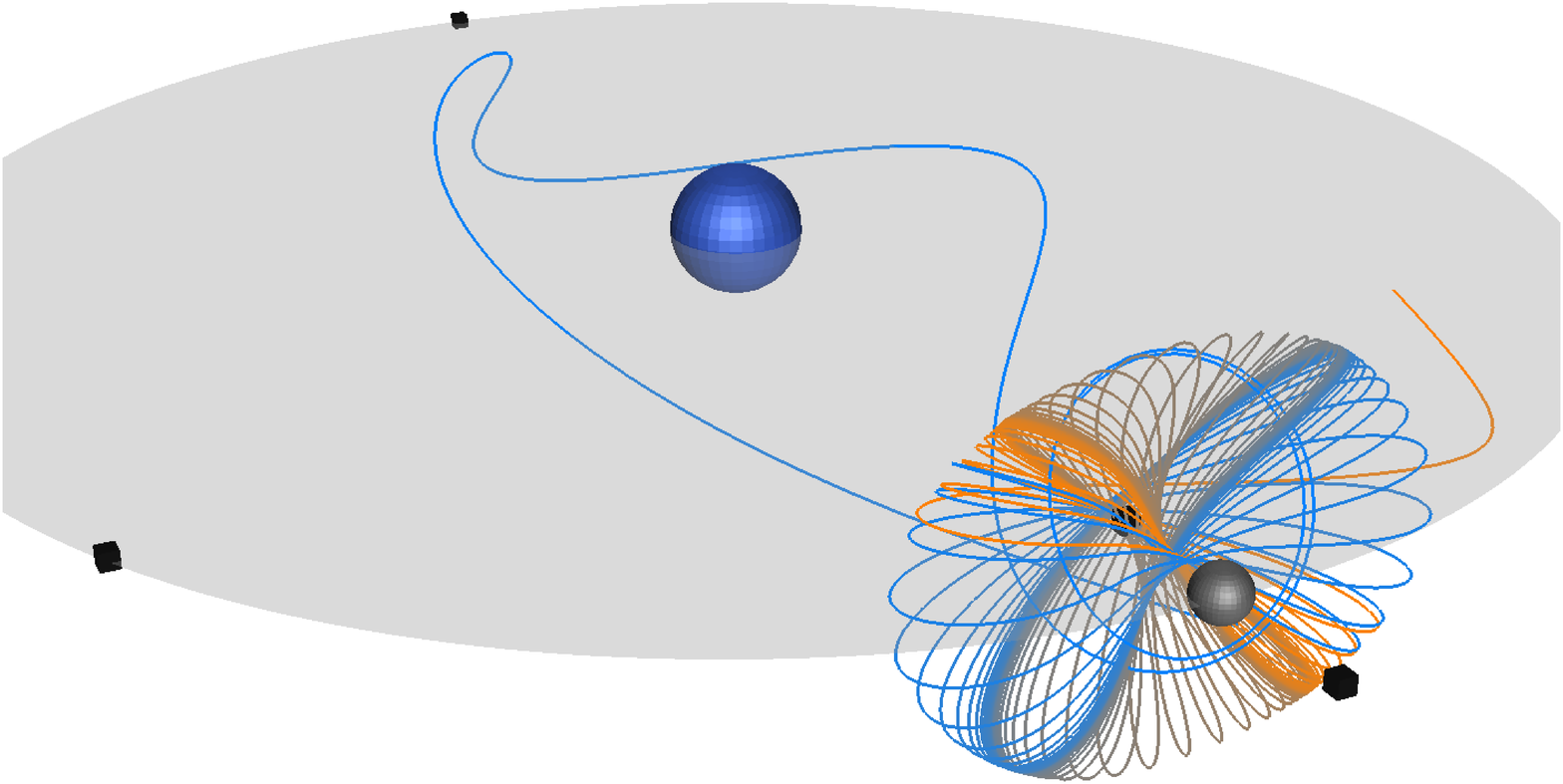}
}
\end{center}
\caption{
(a):
The unstable manifold of the $\H_1$ periodic orbit from Fig.~\ref{fig:manH1},
together with a superimposed longer orbit that was found by continuation using
the $18$-dimensional system in Eq.~\eqref{eq:All}.
(b):
A separate view of the orbit in the unstable manifold. This orbit returns near
the originating $\H_1$ periodic orbit, where it winds around a
quasi-Halo orbit,
while spending significant time near a northern Axial $\A_1$ orbit and its symmetric southern counterpart.
As is clearly visible, and as required by the computational formulation, the orbit
ultimately returns to the plane $\Sigma$, located at $x_{\Sigma}=0.6$.}
\label{fig:H1-to-Axial}
\end{figure}

During the same continuation one encounters the connecting orbit seen in
Fig.~\ref{fig:H1-to-Axial}.
There is still ample evidence of an underlying torus as this figure
shows.
However, on this torus the connecting orbit oscillates between
a northern Axial $\A_1$ orbit and its symmetric southern counterpart,
each time spending a significant number of rotations
very close to each of these two periodic orbits. This suggests a low cost
transfer orbit between an $\H_1$ Halo orbit and the northern or southern
$\A_1$ Axial orbit with the same energy $E=-1.5085$. The connecting
orbit between the northern and southern Axial orbits is generic, as will
be shown in the next section.
Evidence of tori close to heteroclinic connections between two
symmetric Axial families can also be found in Figure 1 of~\cite{Gom-Mon-01}
for an energy value (corresponding to Fig.~\ref{fig:H1-poincare}
discussed later) not far from that of Fig.~\ref{fig:H1-to-Axial}.
In~\cite{Gom-Mon-01} the Axial orbits are denoted by diamond-shaped fixed points.

Fig.~\ref{fig:H1-1Loop}(e) depicts
the change of $\eps$ versus the energy, where the labels correspond
to the values of panels (a)--(d) and Figs.~\ref{fig:H1-to-H1-torus},
\ref{fig:H1-to-Axial}, and~\ref{fig:H1-poincare}.
Note that for most
energy values in this range there exist two connections with different values
of $\eps$.
The continuation terminates at two points on the right hand side of this
diagram, because the Newton-Chord method that AUTO employs no longer
converges there. One of these points corresponds to the special connecting orbit in
Fig.~\ref{fig:H1-to-Axial}. For the other termination point
the unstable real Floquet
multiplier reaches the unit circle, so the two-dimensional unstable
manifold ceases to exist.
Fig.~\ref{fig:H1-1Loop}(f) and the thin-to-thick curve transition for
$\H_1$ in Fig.~\ref{fig:families}(a) show this Floquet multiplier behavior.

\begin{figure}[Htbp]
\begin{center}
\subfloat[]{
\includegraphics[scale=0.198,clip,trim=1mm 13mm 10mm 0mm]{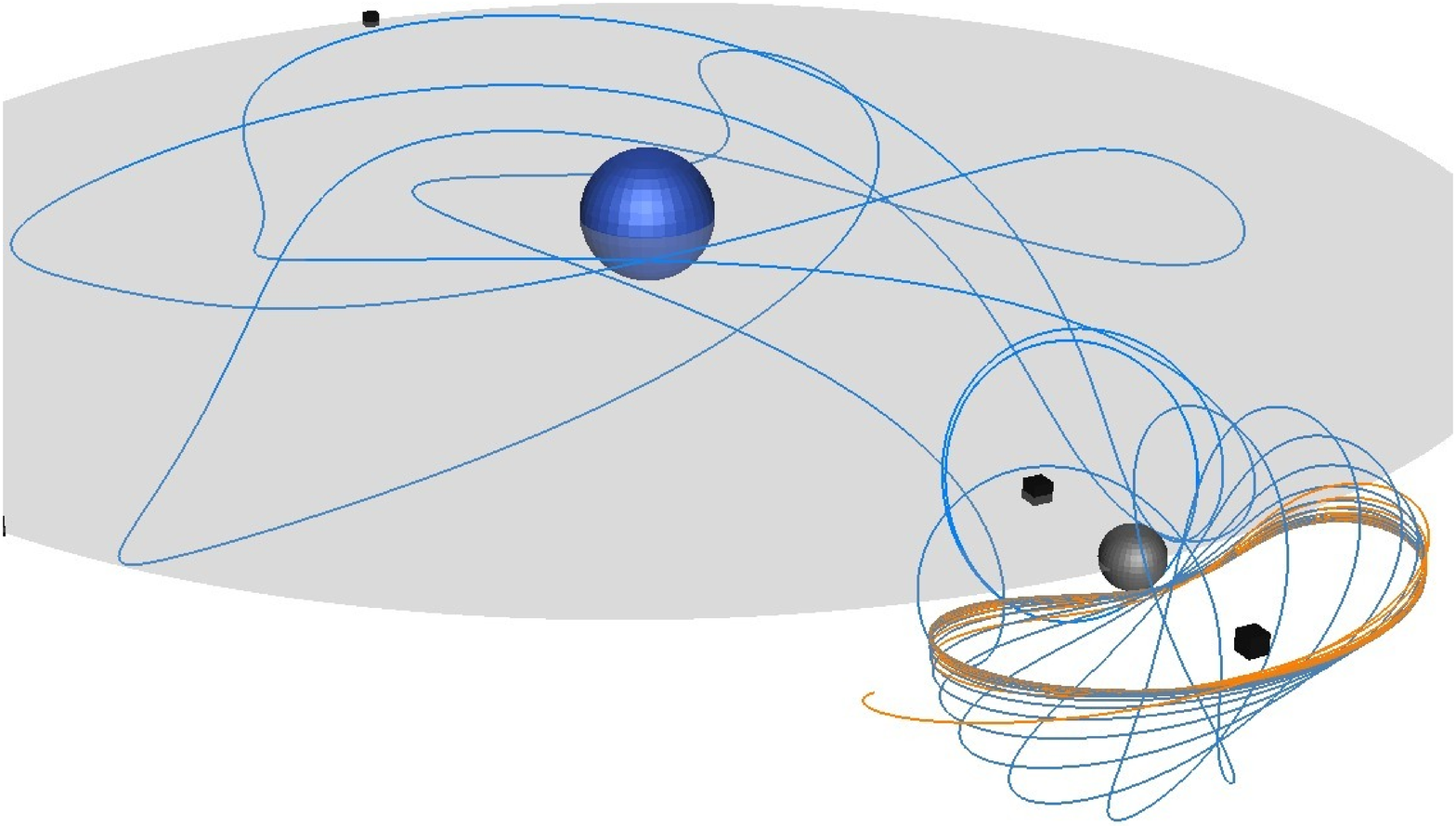}
}
\subfloat[]{
\includegraphics[scale=0.198,clip,trim=1mm 13mm 10mm 0mm]{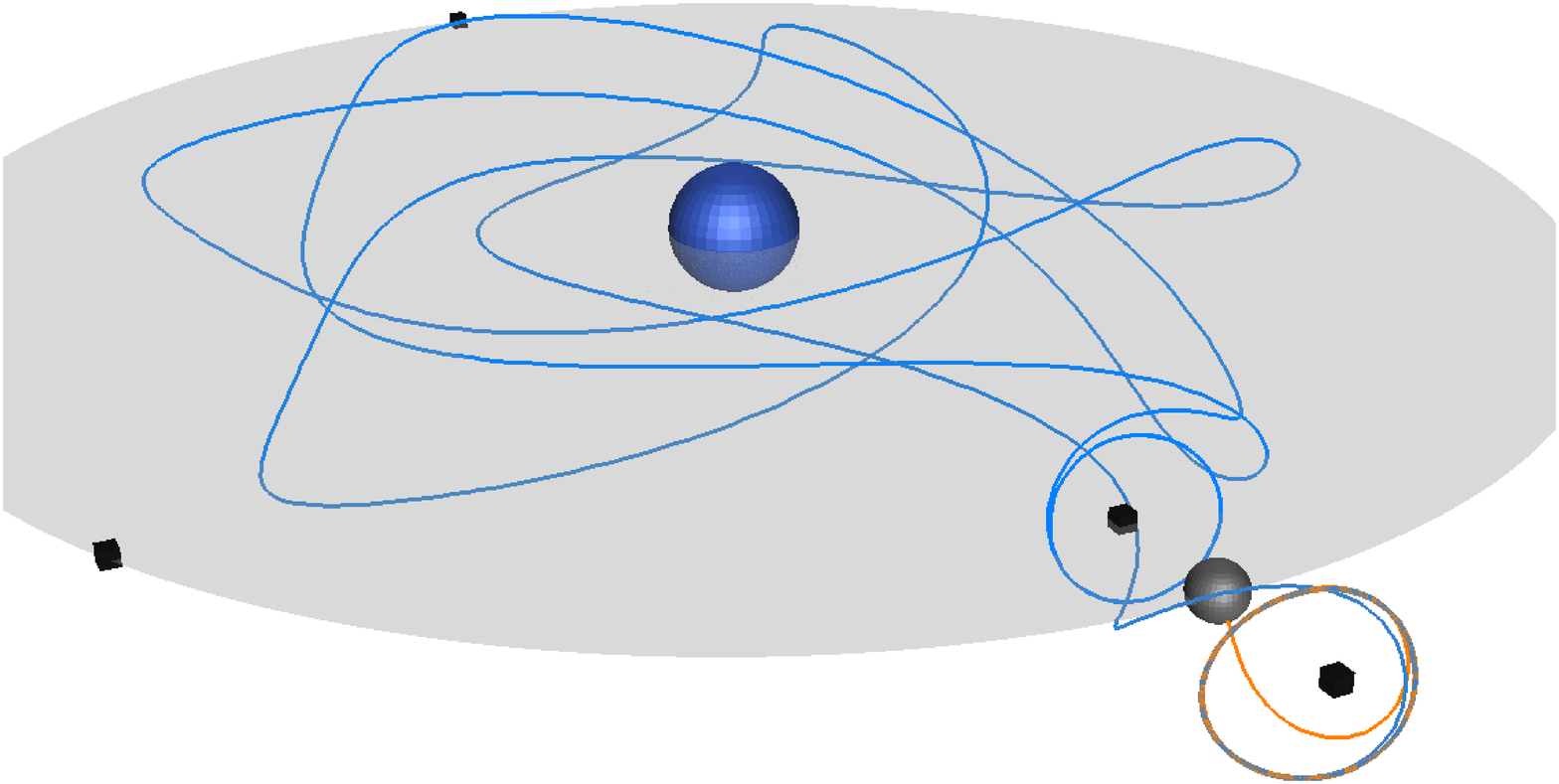}
}

\subfloat[]{\includegraphics[scale=0.198,clip,trim=1mm 13mm 10mm 0mm]{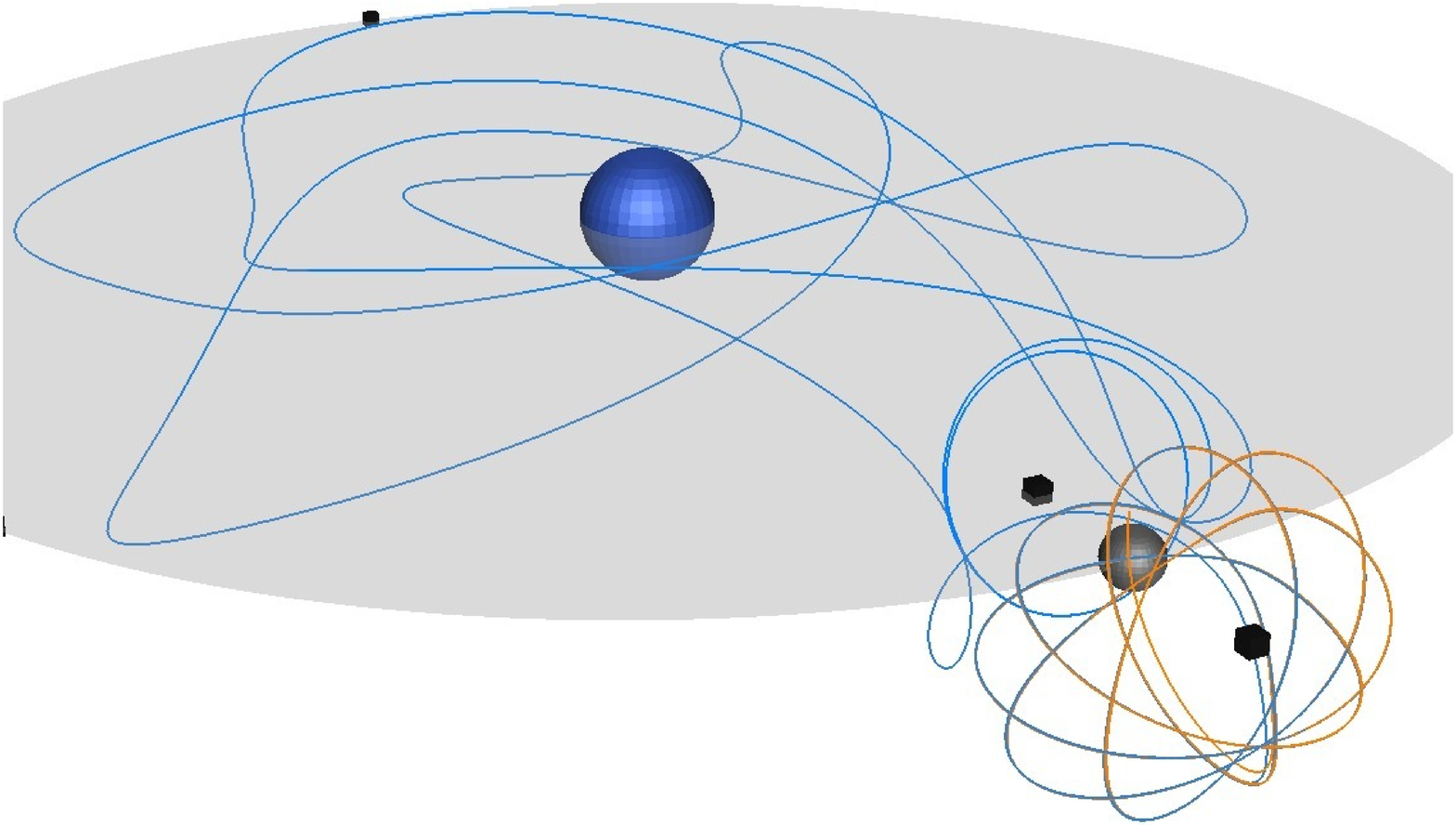}}
\subfloat[]{\includegraphics[scale=0.198,clip,trim=1mm 13mm 10mm 0mm]{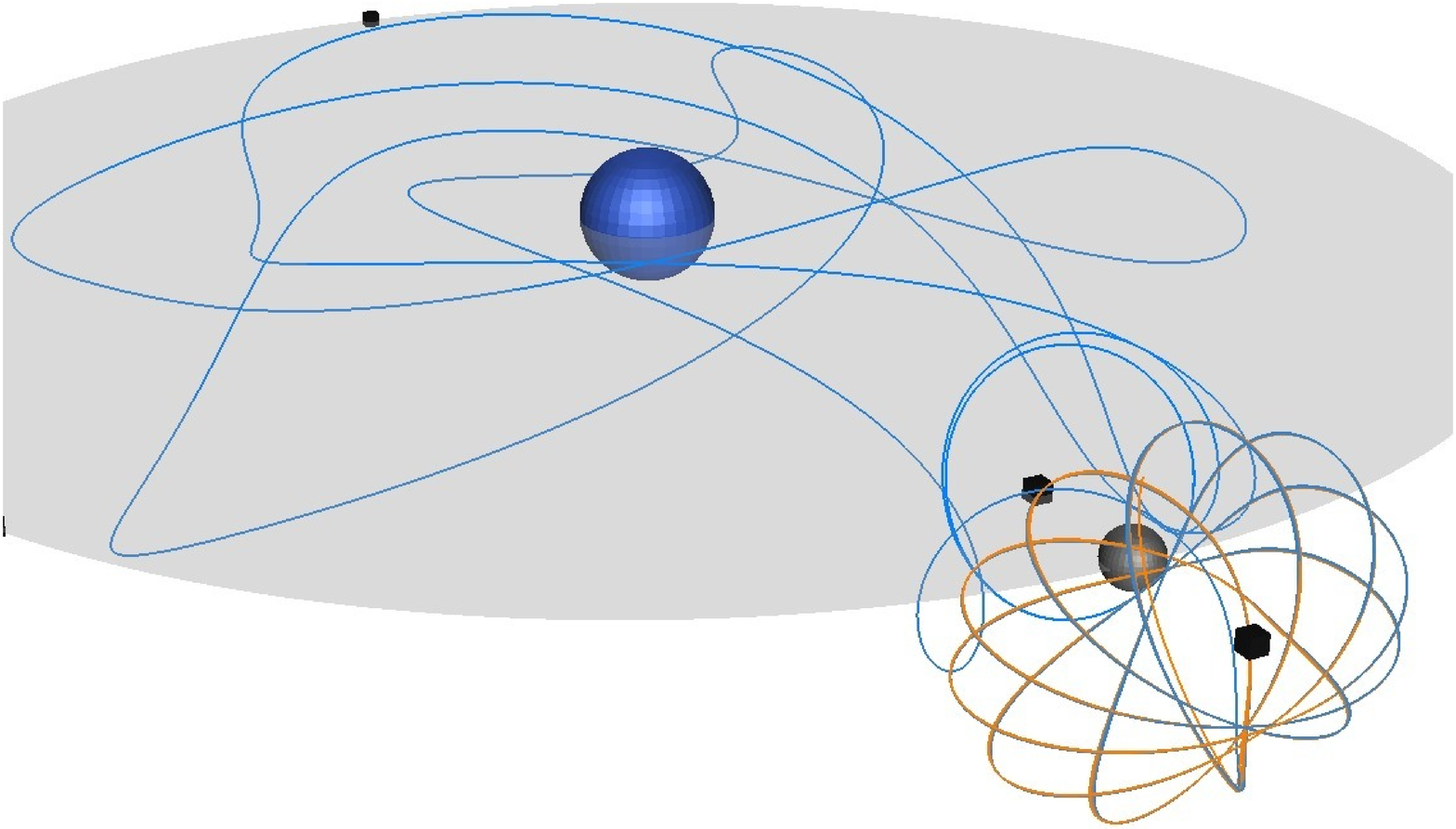}}

\subfloat[]{
\includegraphics{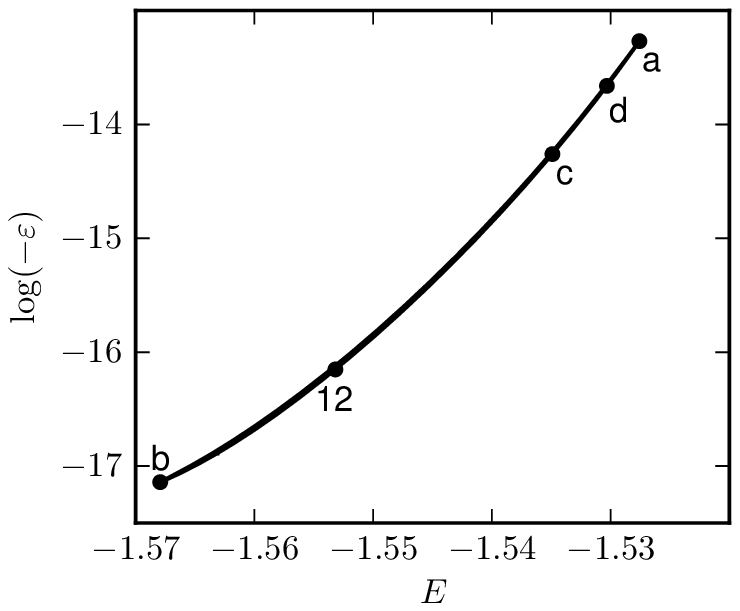}}
\subfloat[]{
\includegraphics{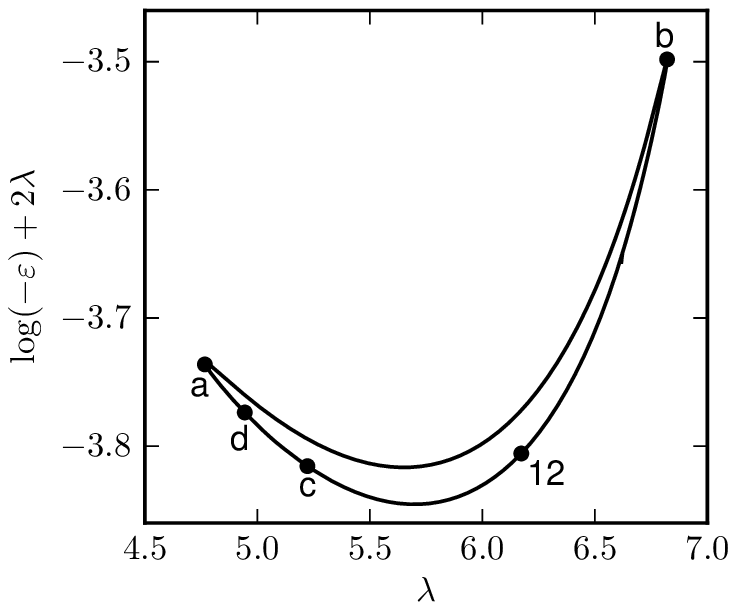}}
\end{center}
\caption{
A family of connecting orbits from $\H_1$ Halo orbits to tori, making
four loops around the Earth.
The terminating plane $\Sigma$ is located at $x_{\Sigma}=1.02$.
(a)--(d): Connecting orbits from
$\H_1$ Halo orbits to (a) a planar $\L_2$ Lyapunov orbit,
(b) approximately an $H_2$ Halo orbit,
(c) a 5:1 resonant orbit near the libration point $\Lib_2$, and
(d) a 6:1 resonant orbit.
(e):
Energy versus $\log(-\eps)$ along the continuation path.
The continuation curve is a loop.
Labels a--d and~\ref{fig:H1-to-H2-torus}
correspond to the values in the panels above and in
Fig.~\ref{fig:H1-to-H2-torus}.
(f):
The logarithm $\lambda$ of the unstable Floquet multiplier of the Halo orbit
versus $\log(-\eps)+2\lambda$ along the continuation path. Here
$\log(-\eps)+2\lambda$ is plotted instead of $\log(-\eps)$ to make the
loop visible.
}
\label{fig:H1-4Loops}
\end{figure}

The apparent connections from an $\H_1$ Halo orbit to an $\L_2$ planar orbit, an $\H_2$ Halo
orbit, and 5:1 and 6:1 resonant orbits shown in
Fig.~\ref{fig:H1-4Loops}(a)--(d)
result from the continuation of the connecting orbit in Fig.~\ref{fig:H1-to-H2-torus},
again using the $18$-dimensional system in Eq.~\eqref{eq:All}.
In this case the connecting orbit makes four loops (instead of just
one as in Fig.~\ref{fig:H1-1Loop}) around the Earth.
Here also, the continuation passes through many resonances. Specifically, the
connecting orbits shown in panels (c) and (d) of Fig.~\ref{fig:H1-4Loops}
approach a period-5 orbit and a period-6 orbit, respectively.

We reiterate that the originating $\H_1$ Halo orbit changes during the continuation,
and with it its period and energy. Consequently the energy of the connecting orbit
as well as the energy of the torus that it approaches change with it, as these energies 
are identical. Similar to panels (e) and (f) of Fig.~\ref{fig:H1-1Loop},
Fig.~\ref{fig:H1-4Loops} show diagrams for the energy and Floquet exponents. Note 
however, that in this case the continuation does not terminate but loops around,
and for every energy value between the two extrema there exist two
connecting orbits. Numerically it was observed that the value of
$\log(-\eps)$ is always close to $-2\lambda$, so to visualize the loop
these two quantities were subtracted in panel (f).

Note that the two extrema in panels (e) and (f) are detected as folds by the 
continuation software, and that these could in turn be followed, for instance by
adding the mass-ratio $\mu$ as a free parameter. The folds themselves do not 
appear to correspond to particularly interesting orbits 
(panels a and b in Fig.~\ref{fig:H1-4Loops} are close to but not at the fold 
locations).

\begin{figure}[Htbp]
\begin{center}
\subfloat[]{
\includegraphics[scale=0.33,clip,trim=40mm 65mm 40mm 50mm]{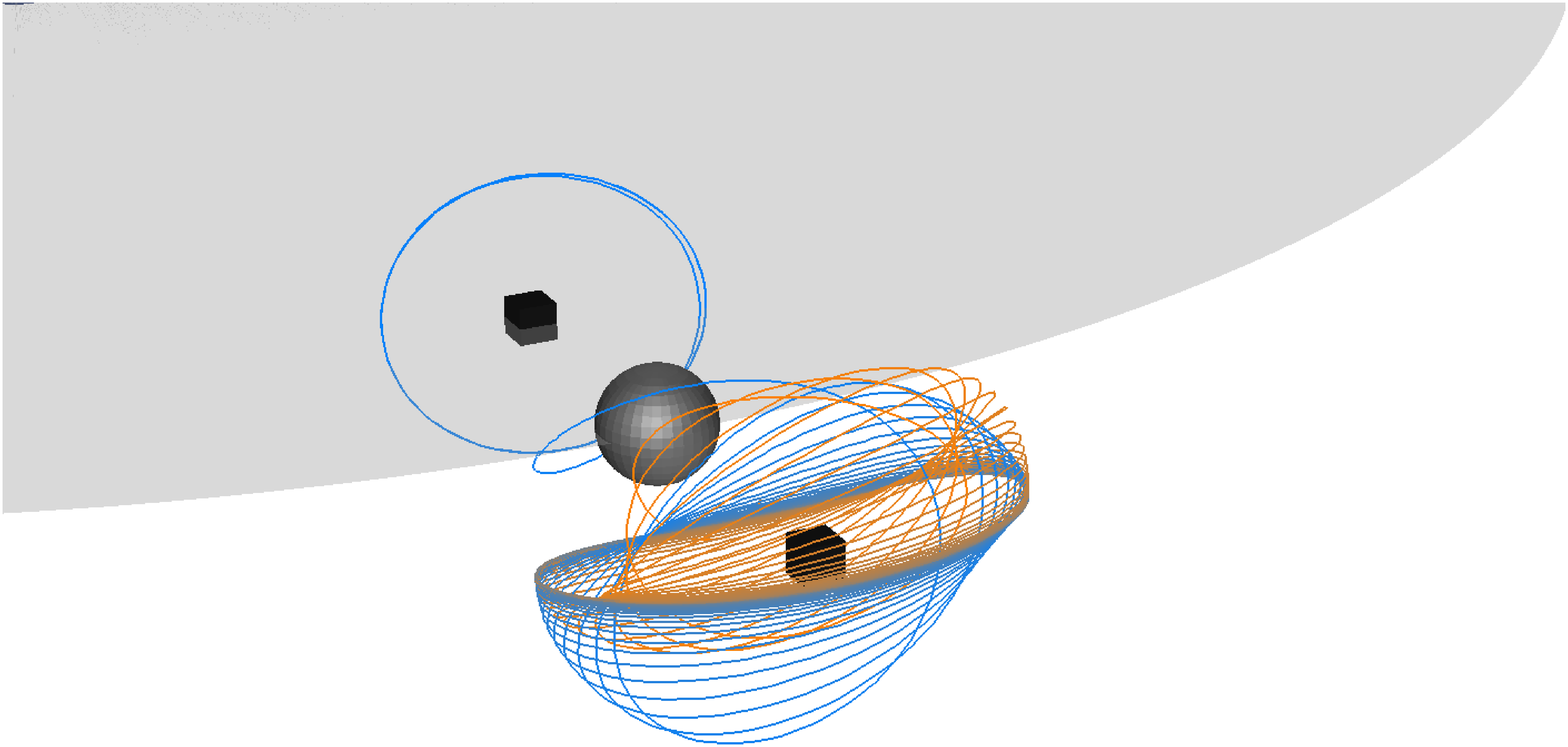}
}
\end{center}
\vskip-0.3truein\noindent
\begin{center}
\subfloat[]{
\includegraphics{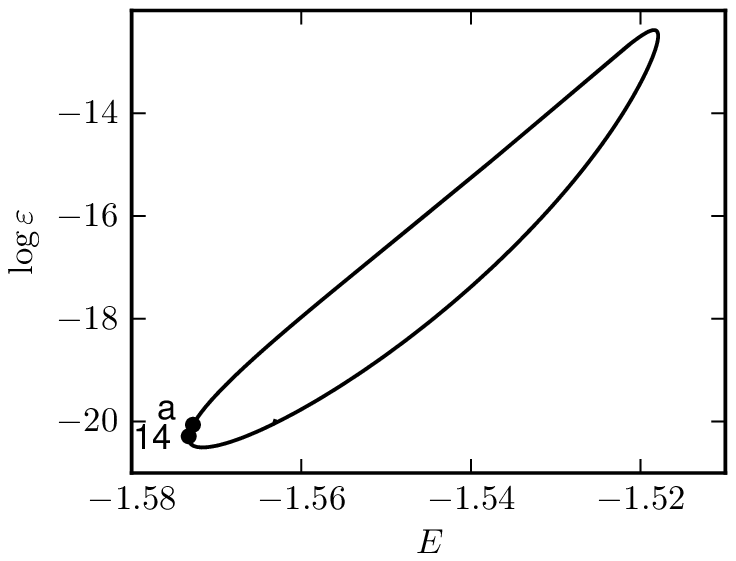}
}
\subfloat[]{
\includegraphics{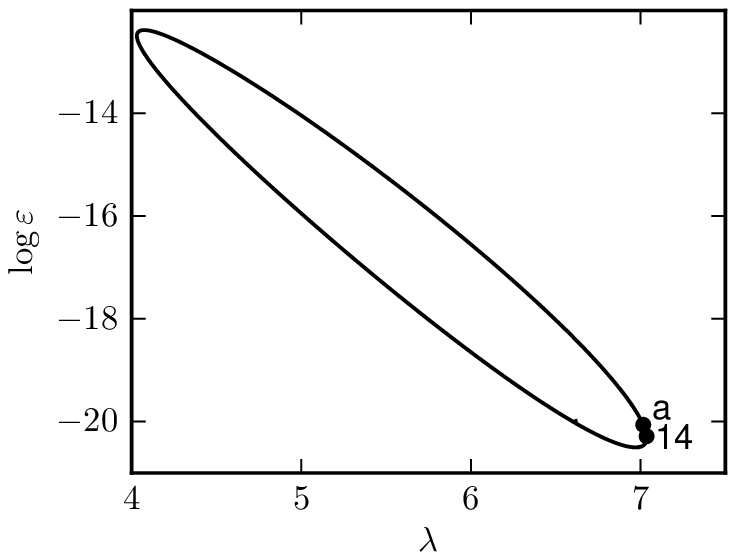}
}
\end{center}
\caption{
(a): An orbit directly connecting
the Moon side of the unstable manifold of the $\H_1$ orbit
with $E=-1.5728$ to an $\L_2$ Lyapunov orbit, obtained by continuation of
the orbit shown in Fig.~\ref{fig:H1-to-H2-Direct}.
(b):
Energy versus $\log(\eps)$ along the continuation path.
The continuation curve is a loop.
Labels a and~\ref{fig:H1-to-H2-Direct}
correspond to the values in the panels above and in
Fig.~\ref{fig:H1-to-H2-Direct}.
(c):
The logarithm of the unstable Floquet multiplier of the Halo orbit
versus $\log(\eps)$ along the continuation path.
}
\label{fig:H1-to-H2-again}
\end{figure}

Similarly, Fig.~\ref{fig:H1-to-H2-again}(a) shows an
interesting orbit from the continuation of the connecting orbit in
Fig.~\ref{fig:H1-to-H2-Direct}, namely a direct connection, \emph{without} looping
around the Earth, from an $\H_1$ orbit to an $\L_2$ Lyapunov orbit.
Fig.~\ref{fig:H1-to-H2-again}(b) and (c) show diagrams
much like those in Fig.~\ref{fig:H1-4Loops}, where the continuation curve
is a loop. However, compared to
Fig.~\ref{fig:H1-4Loops}, here the range of the energy level
and the differences between values of $\log\eps$ are greater.
Also note that the sign of $\eps$ is always positive, rather than negative as before,
as we are now considering the Moon side of the unstable manifold of the $\H_1$ orbits
and so the sign of $\eps$ in Eq.~\eqref{eq:All} is changed.

\section{Existence of Connecting Orbits}
\label{sec:discussion}

As seen in the preceding sections, boundary value formulations
provide a powerful tool to compute unstable manifolds and
connecting orbits in the CR3BP. Our exploration shows only
the tip of the iceberg of a wealth of interesting orbits. This section
discusses how the connecting orbits that were found relate to the existing
literature.
As mentioned in Sect.~\ref{sec:introduction}, the periodic orbits and tori
themselves, as well as homoclinic and heteroclinic
orbits connecting $\L_1$ and $\L_2$ periodic orbits have been
studied extensively. However, for the spatial case, to the best of our
knowledge, connecting orbits from $\H_1$ Halo to quasi-$\H_1$ or quasi-$\H_2$ have not
been explicitly found before. We must mention that connections from
quasi-$\H_1$ to quasi-$\H_2$ orbits can be found in 
\cite{Gom-Koo-Lo-Mar-Mas-Ros-04}.

\begin{table}[Htbp]
\begin{center}
{\normalsize
\begin{tabular}{|l|l|l|l|l|}
\hline
Figure & Energy & Type & Period & Non-trivial Floquet exponents \\
\hline
\ref{fig:manV1},\ref{fig:man+orbit-V1},\ref{fig:orbit-V1},\ref{fig:shooting}
& $-1.5164$ & $\V_1$ & 3.7700 &
$\pm 6.4948$, $\pm 0.077175 \cdot 2\pi i$ \\
\hline
\ref{fig:H1-to-H2-torus} & $-1.5532$ & $\H_1$ & $2.7873$ &
$\pm 6.1730$, $\pm 0.21859\cdot 2\pi i$\\
& & quasi-$\H_2$ & $3.2753$ & $\pm 5.9348$, $\pm 0.17737 \cdot 2\pi i$\\
\hline
\ref{fig:H1-to-H1-torus} & $-1.5631$ & $\H_1$ & $2.7821$&
$\pm 6.6179$, $\pm 0.16367\cdot 2\pi i$\\
\hline
\ref{fig:H1-to-H2-Direct} & $-1.5733$ & $\H_1$ & $2.7716$
& $\pm 7.0356$, $\pm 0.11217 \cdot 2\pi i$ \\
& & quasi-$\H_2$ & $3.3805$ & $\pm 6.7839$, $\pm 0.068079 \cdot 2\pi i$ \\
\hline
\ref{fig:H1-1Loop}(a) & $-1.5552$ & quasi-$\H_1$ (5:1 res) &
$2.7868$ & $\pm 6.2673$, $\pm 0.20690\cdot 2\pi i$\\
\hline
\ref{fig:H1-1Loop}(b) & $-1.5716$ & $\H_1$ & $2.7736$ &
$\pm 6.9682$, $\pm 0.12050\cdot 2\pi i$ \\
\hline
\ref{fig:H1-1Loop}(c) & $-1.5754$ & $\H_1$ & $2.7690$ &
$\pm 7.1183, \pm 0.10189\cdot 2\pi i$  \\
& & $\L_1$ & $2.8982$ & $\pm 0.29251$, $\pm 7.4268$\\
\hline
\ref{fig:H1-1Loop}(d) & $-1.5617$ & $\H_1$ & $2.7832$ &
$\pm 6.5556$ ,$\pm 0.17135\cdot 2\pi i$\\
\hline
\ref{fig:manH1},\ref{fig:H1-to-Axial} & $-1.5085$ & $\H_1$ & $2.5152$ &
$\pm 2.8541, \pm 0.38928 \cdot 2\pi i$ \\
& & $\A_1$ & $4.0117$ & $\pm 0.43035$, $\pm 6.1278$ \\
\hline
\ref{fig:H1-4Loops}(a) & $-1.5276$ & $\H_1$ & $2.7472$ &
$\pm 4.7666, \pm 0.39460 \cdot 2\pi i$\\
& & $\L_2$ & $3.9550$ & $\pm 5.9229$, $\pm 0.52353$ \\
\hline
\ref{fig:H1-4Loops}(b) & $-1.5679$ & $\H_1$ & $2.7776$ &
$\pm 6.8207$, $\pm 0.13870 \cdot 2\pi i$ \\
& & $\H_2$ & $3.3567$ & $\pm 6.5802$, $\pm 0.095719 \cdot 2\pi i$\\
\hline
\ref{fig:H1-4Loops}(c) & $-1.5349$ & $\H_1$ &
$2.7715$ & $\pm 5.2225, \pm 0.33730 \cdot 2\pi i$ \\
& & quasi-$\H_2$ (5:1 res) & $3.1170$ & $\pm 4.9006$, $\pm 0.30348 \cdot 2 \pi i$ \\
\hline
\ref{fig:H1-4Loops}(d) & $-1.5303$ & $\H_1$ &
$2.7581$ & $\pm 4.9440, \pm 0.23394 \cdot 2\pi i$\\
& & quasi-$\H_2$ (6:1 res) & $3.0568$ & $\pm 4.5639$, $\pm 0.34345 \cdot 2 \pi i$ \\
\hline
\ref{fig:H1-to-H2-again}(a) & $-1.5728$ & $\H_1$ &
$2.7722$ & $\pm 7.0158$,$\pm 0.11462\cdot 2\pi i$ \\
& & quasi-$\H_2$ & $3.3784$ & $\pm 6.7651$, $\pm 0.070727 \cdot 2\pi i$\\
\hline
\ref{fig:H1-poincare} & $-1.5130$ & $\H_1$ & $2.6176$ & $\pm
3.5327$,$\pm 0.45736 \cdot 2 \pi i$ \\
\hline
\end{tabular}
}
\caption{Numerically computed values for all figures. The first line in
  each box contains the data for the starting periodic orbit $\u(t)$.
  If the orbit $\r(t)$ connects the periodic orbit to a different
  periodic orbit, or to a torus that surrounds a different periodic
  orbit (a quasi-periodic orbit), then the second row contains data
  for this second periodic orbit. There are always two trivial Floquet
  exponents equal to zero.}
\label{tab:Erg-Typ-Per-Flq}
\end{center}
\end{table}

In all three continuations in Sect.~\ref{sec:connecting-orbits-examples}
a collection of interesting objects was seen, including apparent heteroclinic 
connections
between Halo orbits and resonant orbits, {\it i.e.}, seemingly heteroclinic
connections between Halo orbits and $n$-periodic orbits, where
the complex Floquet multipliers of the Halo orbit near the
$n$-periodic orbit are close to $e^{2\pi i/n}$.
Numerical data corresponding to figures in this paper are given to 5 significant 
digits in Table~\ref{tab:Erg-Typ-Per-Flq}.

The existence of certain connecting orbits can be explained by counting dimensions.
For Fig.~\ref{fig:H1-to-Axial},
the multipliers of the symmetric Axial orbits give rise to
three-dimensional stable and unstable manifolds, so the connection
between the two northern and southern $\A_1$ Axial orbits seen
in Figure~\ref{fig:H1-to-Axial} is generic for the CR3BP posed in $\mathbb{R}^6$.
The connection from the
$\H_1$ Halo orbit to the $\A_1$ Axial orbit is codimension-one, since the $\H_1$  Halo orbit has a two-dimensional
unstable manifold, and the stable manifold of the $\A_1$ Axial orbit is three-dimensional.
Similarly the planar $\L_1$ and
$\L_2$ Lyapunov orbits for the energy values that we observe
in Figs.~\ref{fig:H1-1Loop}(c), \ref{fig:H1-4Loops}(a), and
\ref{fig:H1-to-H2-again}(a) have
three-dimensional stable and unstable manifolds, so connections
from $\H_1$ Halo Orbits to those orbits are also codimension-one. 
By standard bifurcation theory, codimension-one connecting orbits can be 
expected to arise as one parameter is varied, for specific values of that 
parameter. This is indeed what we observed in these cases, where we found
such connecting orbits for specific values of the energy $E$.

\subsection{Connecting Orbits from a Halo orbit to a quasi-Halo orbit}

In contrast to the connecting orbits considered above,
the connections from the $\H_1$ Halo orbit to the same or other Halo orbits
cannot simply be explained by counting dimensions. The $\H_1$ Halo orbit
has a two-dimensional unstable manifold and a two-dimensional stable manifold,
so such connections would be codimension-two and would not be expected
along a family (in the sense of this paper) of Halo orbits.
The apparent Halo orbits in the Halo to Halo connecting orbits
may in fact be very thin quasi-Halo orbits
(tori). However, even if this distinction were mathematically important,
it would resumably be of little practical importance in space mission design.

\begin{figure}[htbp]
\subfloat[]{
\includegraphics[scale=0.29,clip,trim=90mm 16mm 50mm 5mm]{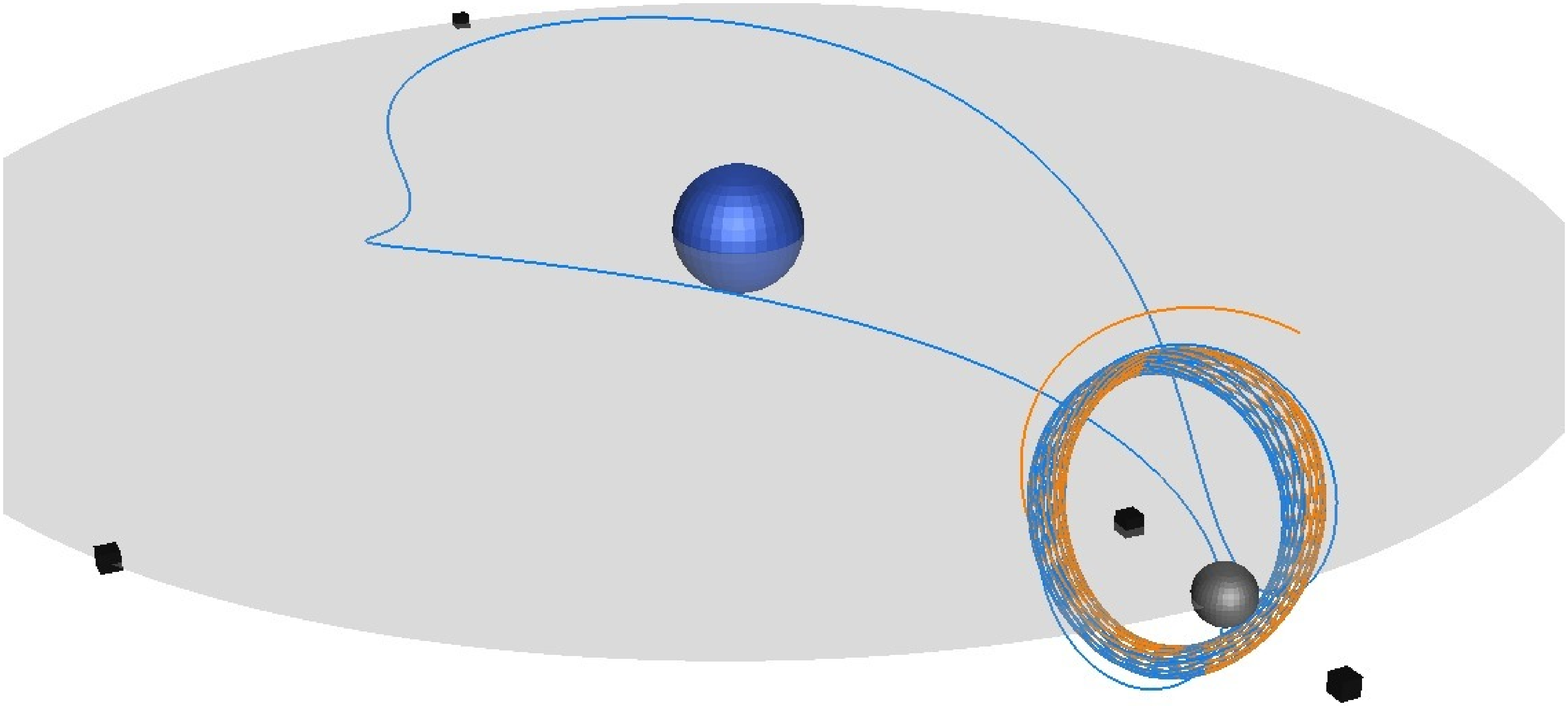}
}
\subfloat[]{
\includegraphics[scale=0.29,clip,trim=106mm 16mm 12mm 5mm]{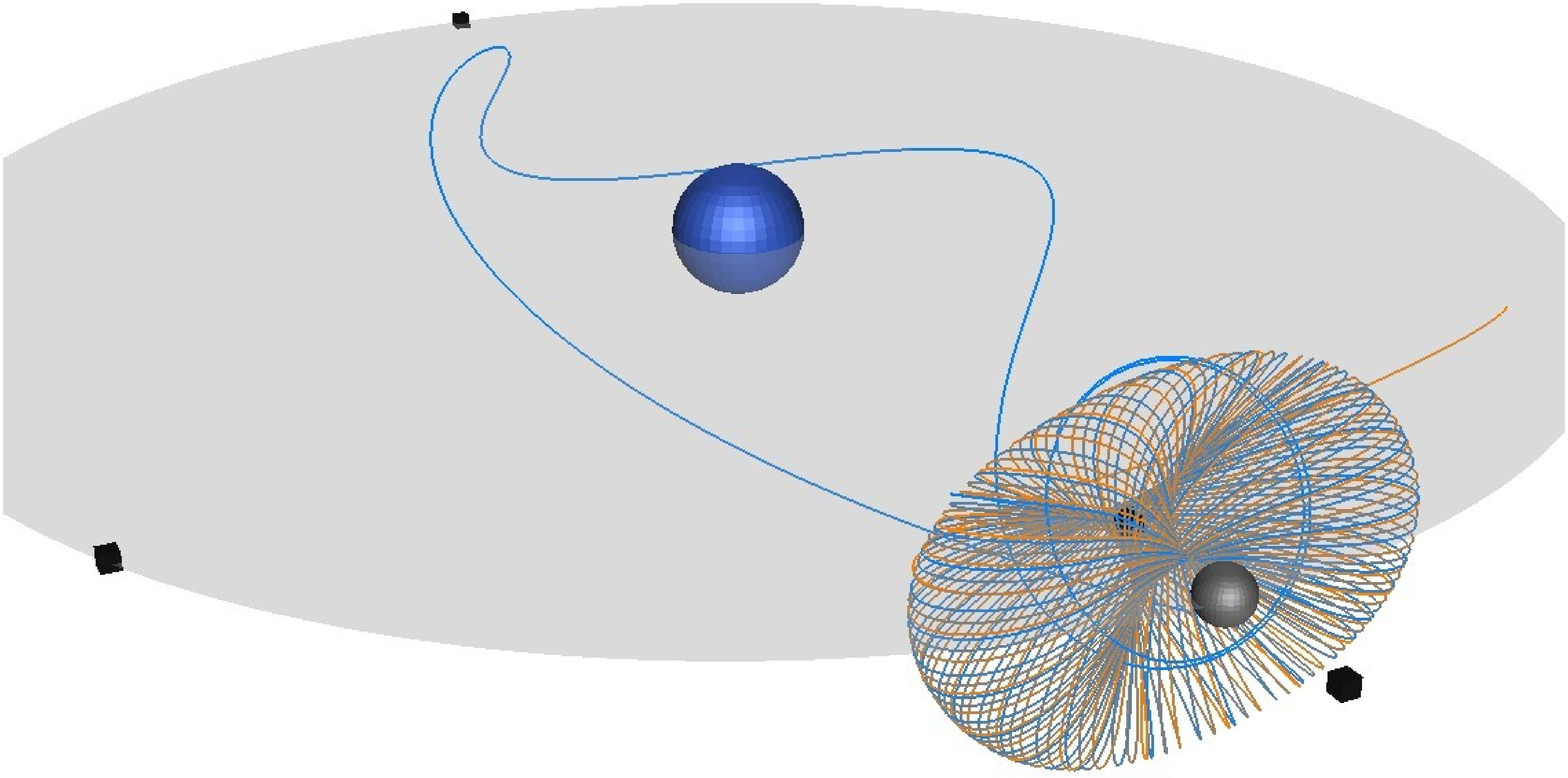}
}
\vskip-0.1truein
\subfloat[]{
\includegraphics{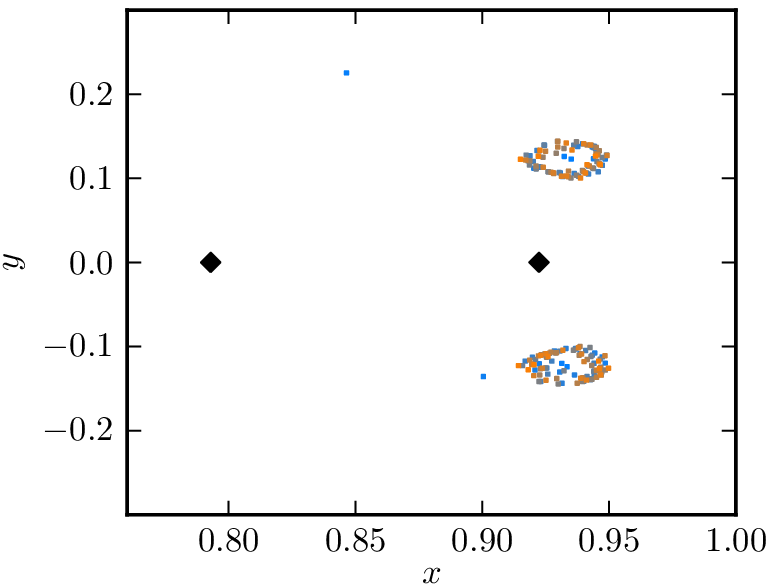}
}
\subfloat[]{
\includegraphics{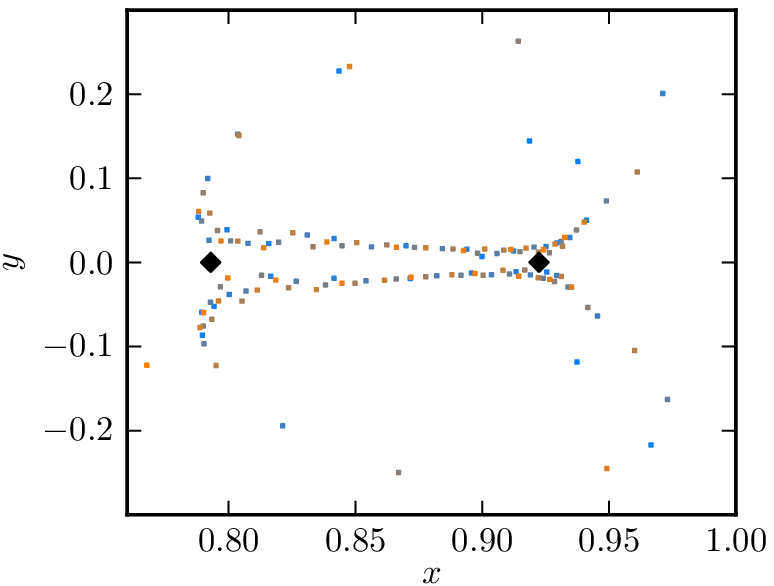}
}
\caption{(a) and (b): Two orbits for the family displayed in
Fig.~\ref{fig:H1-1Loop} where $E=-1.5130$. (c) and (d):
Their intersections with the plane $z=0$.
The energy value corresponds to Figure~1 in~\citet{Gom-Mon-01}
and Figure~3 in~\citet{Ale-Gom-Mas-09},
where the energy is denoted as $h=E+\mu(1-\mu)/2=-1.507$.
The two diamonds correspond to intersections of the $\A_1$ Axial orbits,
as in \citet{Gom-Mon-01}.}
\label{fig:H1-poincare}
\end{figure}

We justify the presence of
the connecting orbits from $\H_1$ Halo orbits to quasi-Halo
orbits by recalling the detailed analysis of
center manifolds around collinear libration points in
\citet{Jor-Mas-99}, \citet{Gom-Koo-Lo-Mar-Mas-Ros-04}, and \citet{marsden}.
For a prescribed energy level, the thin solid curves in
Fig.~\ref{fig:families} show that
there exists a Halo orbit near $\Lib_i$, for $i=1,2$, that has two
complex Floquet multipliers on the unit circle and different from one.
There are two more real Floquet multipliers, one
of them greater than one and the other less than one. The Floquet multipliers
and corresponding energy levels for the Halo orbits in Figs.~\ref{fig:manH1}
and~\ref{fig:H1-to-H2-torus} to~\ref{fig:H1-poincare} are shown
in Table~\ref{tab:Erg-Typ-Per-Flq}.

The unstable manifold of such a Halo orbit is a two-dimensional
surface, see Fig.~\ref{fig:manH1}.
On the other hand, close to either of the two libration points $\Lib_j$,
for $j=1,2$, there is
a four-dimensional center manifold, which exists
due to the fact that $\Lib_j$ has four eigenvalues on the unit circle and two
more on the real line excluding one.

On the energy surface the center manifold of the co-linear libration
points is a little different.
In \citet{marsden} the authors
verify that for fixed energy, the center manifold of $\Lib_j$, for $j =1,2$, is a
normally hyperbolic invariant manifold (NHIM) corresponding
to a normally hyperbolic
three-sphere that is invariant for the linearized system.
Thus, in a small neighborhood of  $\Lib_j$ the center
manifold becomes a deformed three-sphere for the nonlinear system.
Moreover, it is a three-dimensional hyperbolic
invariant manifold in a five-dimensional energy surface with one stable
direction. Therefore the stable manifold of the center manifold is
four-dimensional.

In this setting, dimension counting
gives a heuristic argument that motivates
why connecting orbits from
Halo orbits to torus-like objects appear
quite naturally in our problem, as seen in Figs.~\ref{fig:H1-to-H2-torus},
\ref{fig:H1-to-H1-torus}, \ref{fig:H1-to-H2-Direct}, and~\ref{fig:H1-poincare}.
We have mentioned that the center
manifold of a libration point $\Lib_i$ restricted to
the energy surface
is a three-dimensional NHIM,
so we can do a dimension counting analysis similar to those
above.
The main observation is that for a fixed value of the
energy, these quasi-Halo
orbits are lower dimensional whiskered quasi-periodic
solutions (see \citet{Fon-Lla-Sir-09}) that lie inside the three-dimensional
center manifolds of the collinear libration points.

Now, the Halo orbit is a one-dimensional normally hyperbolic object with one
unstable direction. The Halo orbit itself also belongs to the center manifold
of $\Lib_i$ and its unstable manifold is a two-dimensional object in the energy surface.
Finally, by dimension counting, we notice that a transversal intersection
between the unstable manifold of the Halo orbit of $\Lib_j$ and the stable manifold of
the center manifold of $\Lib_j$ is a one-dimensional object.

Several authors \citep{Jor-Mas-99,Mas-05,Ale-Gom-Mas-09}
have explored the center manifolds of collinear libration points
by means of Lindstedt-Poincar\'e expansions
of the Hamiltonian restricted to each center manifold.
These computations suggest that for the energy levels under
consideration in this paper there exist large regions in the
center manifold that exhibit quasi-periodic motions.
The $z=0$ sections of the center manifolds 
reveal that most of these quasi-periodic obits are two-dimensional
invariant tori in the center manifold
(see for instance Figure 3 in \citet{Ale-Gom-Mas-09}).
These invariant tori are also normally hyperbolic and their corresponding stable
manifolds are two-dimensional objects inside the four-dimensional stable manifold
of the center manifold. These codimension-one stable manifolds separate the
four-dimensional stable manifold into regions.

A trajectory in the intersection of the unstable manifold of the Halo orbit and the
stable manifold of the center manifold has a positive probability
of being either on the stable manifold of
an invariant torus or in a region bounded by stable manifolds of
invariant tori. It is even possible that the trajectory belongs to the
stable manifold of a periodic orbit ({\it e.g.}, a Halo orbit) in the center manifold 
of $\Lib_j$. However, this is unlikely since the stable manifold of a periodic orbit 
is a two-dimensional object inside a four-dimensional stable manifold. Thus, it is more 
likely that a trajectory in the intersection of the unstable manifold of a Halo orbit 
with the stable manifold of the three-dimensional center manifold of a libration point 
passes close to a two-dimensional hyperbolic quasi-periodic orbit in the center manifold.

Fig.~\ref{fig:H1-poincare}(a) and (b) show two connections
from a Halo orbit to a quasi-Halo orbit nearby, from the continuation
displayed in Fig.~\ref{fig:H1-1Loop}, where the energy value
$E=-1.5130$ matches the value $h=-1.507$ used for Figure 1 in
\citet{Gom-Mon-01} and Figure 3 in \citet{Ale-Gom-Mas-09},
where $h=E+\mu(1+\mu)/2$. In Fig.~\ref{fig:H1-poincare}(c) and (d)
we show the intersections of these trajectories with the plane $z=0$.
If we compare these intersections to the corresponding ones in
\citet{Gom-Mon-01} and \citet{Ale-Gom-Mas-09},
we discover that the trajectory seems to approach one of these
quasi-Halo orbits and then drift away from it. We conclude that
the trajectory computed by our methods is shadowing a connecting
orbit in the intersection of the unstable manifold of the Halo orbit
and the stable manifold of the
corresponding normally hyperbolic lower dimensional torus.
\section{Numerical Aspects}
\label{sec:numerical-aspects}
In our computations we used the continuation and bifurcation software
AUTO \citep{AUTO81,AUTO} for computing families of periodic orbits, associated
Floquet eigenfunctions, unstable manifolds, and connecting orbits.
Orbits are continued in AUTO as solutions to a suitable boundary value problem
(BVP), as described in this paper.
To compute approximate solutions of BVPs, AUTO uses the method of Gauss-Legendre
collocation with piecewise polynomials on adaptive meshes
\citep{deBoorSwartz,AscherChristiansenRussell}. 
These calculations can be fast; for example the {\it entire} Vertical family of 
periodic orbits $\V_1$ can be computed in 0.12 seconds on a dual core laptop, using 
as few as $20$ mesh intervals with 4 Gauss collocation points each, and using $48$ 
continuation steps. In our actual computations we use more mesh intervals to ensure 
high accuracy, {\it e.g.}, with $100$ mesh intervals the family $\V_1$ can be computed 
in 0.25 seconds. We also often use more continuation steps, for example for the
purpose of generating data for computer animations.

Two-dimensional unstable manifolds computed as a solution family 
by continuation, as done in Sect.~\ref{sec:manifold-algorithm}, require more mesh 
intervals and continuation steps when the orbits in the family wind many times
around a torus. For example, a connection from the Halo orbit $\H_1$ to a torus 
near $\H_2$, winding around the torus 12 times, can be reached easily in 15 seconds, 
using 200 mesh intervals and 1000 continuation steps. However, we also
computed such connections with a much higher number of windings, which requires a 
correspondingly higher number of mesh intervals and continuation steps. 
In fact, we have used up to $1500$ mesh intervals in such cases. Furthermore, to 
ascertain the correctness of our results, we have computed these manifolds with 
various choices of the number of mesh intervals.

Similarly, the continuation of solutions to the $18$-dimensional system to follow 
periodic-orbit-to-torus connections, as given in Sect.~\ref{sec:connecting-orbits} 
and used in Sect.~\ref{sec:connecting-orbits-examples}, requires a correspondingly 
high number of mesh intervals. Since the dimension of the system is then $18$, 
{\it i.e.}, three times the dimension of the systems used for continuing periodic 
orbits and unstable manifolds, the computation time would increase by a factor $3^3$, 
that is $27$. However, this is reduced since AUTO takes into account the zero 
structure of the submatrices of the full Jacobian matrix that correspond to 
individual mesh intervals. On the other hand, the connecting orbit requires 
significantly more mesh intervals then the base periodic orbit and its Floquet 
eigenfunction, and AUTO does not take advantage of this. Thus the continuation 
of the coupled system (periodic orbit, Floquet eigenfunction, connecting orbit) 
in Sect.~\ref{sec:connecting-orbits} can take significant computer time, also 
because such continuation with varying energy requires many continuation steps. 
In extreme cases the calculations have taken up to $6$ hours computer
time.

\begin{acknowledgements}
The work of EJD, ARH, and BEO was supported by NSERC (Canada) Discovery Grants.
RCC acknowledges support by a FQRNT PBEEE (Qu\'ebec) award, and ALR is
supported by a Conacyt (M\'exico) scholarship.
\end{acknowledgements}

\end{document}